\documentclass[10pt,reqno]{amsart}
\usepackage{amsfonts,amsmath,latexsym,verbatim,amscd,mathrsfs,color,array}
\usepackage[colorlinks=true,citecolor=blue]{hyperref}

\usepackage{appendix}
\usepackage{amsmath,amssymb,amsthm,mathrsfs,amsfonts,graphicx,color}
\usepackage{amsthm} 
\usepackage[nameinlink]{cleveref} 
\usepackage{amsthm}

\usepackage{amsmath}

\usepackage{amsfonts,amsmath,latexsym, verbatim,amscd,mathrsfs, color,array}
\usepackage[colorlinks=true,citecolor=blue]{hyperref}

\usepackage{appendix}
\usepackage{amsmath,amssymb,amsthm,amsfonts,graphicx, color}

\newcommand{\cuad}{{\sqcap\kern-.68em\sqcup}}

\newcommand{\be}{\begin{equation}}
	\newcommand{\ee}{\end{equation}}

\renewcommand{\div}{\mathop{\rm div}}

\newtheorem{definition}{Definition}
\newtheorem{lemma}{Lemma}[section]
\newtheorem{prop}{Proposition}[section]
\newtheorem{theorem}{Theorem}[section]
\newtheorem{corollary}{Corollary}[section]
\newtheorem{remark}{Remark}[section]
\newcommand{\bremark}{\begin{remark} \em}
	\newcommand{\eremark}{\end{remark} }

\numberwithin{equation}{section}

\def\ep{\varepsilon}

\title[Stability of 3D Capillary-Gravity Waves]{Stability of Solitary Capillary-Gravity Water Waves in Three Dimensions}

\author[C. Gui]{Changfeng Gui}
\address{\noindent  Changfeng Gui, Department of Mathematics, Faculty of Science and Technology, University of Macau, Taipa, Macao SAR, China.}
\email{changfenggui@um.edu.mo}

\author[S. Lai]{Shanfa Lai}
\address{\noindent  Shanfa Lai, Department of Mathematics, Faculty of Science and Technology, University of Macau, Taipa, Macao SAR, China.
}  \email{laishanfa@amss.ac.cn}

\author[Y. Liu]{Yong Liu}
\address{\noindent  Yong Liu, School of Mathematics and Statistics, Beijing Technology and Business University, Beijing, China.}
\email{yliumath@btbu.edu.cn}

\author[J. Wei]{Juncheng Wei}
\address{\noindent  Juncheng Wei, Department of Mathematics, The Chinese University of Hong Kong, Shatin, Hong Kong.
}  \email{wei@math.cuhk.edu.hk}

\author[W. Yang]{Wen Yang}
\address{\noindent  Wen Yang, Department of Mathematics, Faculty of Science and Technology, University of Macau, Taipa, Macao SAR, China.
}  \email{wenyang@um.edu.mo}

\date{\today}

\begin{document}

	\begin{abstract}
		This paper establishes the conditional orbital stability of fully localized solitary waves for the three-dimensional capillary-gravity water wave problem in finite depth under strong surface tension. The waves, constructed via a non-variational Lyapunov–Schmidt reduction in \cite{GLLWY}, are not energy minimizers and thus require a direct stability analysis. We adapt the Grillakis–Shatah–Strauss framework within Mielke's approach to handle the mismatch between well-posedness and energy spaces. The proof relies on spectral analysis of the linearized dynamics and careful treatment of the Hamiltonian structure defined by the energy and momentum functionals.
	\end{abstract}

	\maketitle

	\section{Introduction}
	
	In this article, we investigate the stability of traveling wave solutions for capillary-gravity water waves in a three-dimensional fluid domain of finite depth. The physical system is governed by the incompressible Euler equations,
	\begin{equation}
		\label{eq:euler1}
		\frac{\partial \mathbf{u}}{\partial t} + (\mathbf{u} \cdot \nabla) \mathbf{u} + \nabla P = g \mathbf{e}_3,
	\end{equation}
	on the evolving fluid domain
	\begin{equation}
		\label{eq:euler2}
		\Omega(t)=\{(x_1, x_2, x_3) \in \mathbb{R}^3 : -1< x_3< \eta(t, x_1, x_2)\}.
	\end{equation}
	Here, the water density is normalized to unity, $\mathbf{u}(t, \cdot): \Omega(t) \to \mathbb{R}^3$ denotes the velocity field, $P(t, \cdot): \Omega(t) \to \mathbb{R}$ is the pressure, $g > 0$ is the gravitational constant, and $\mathbf{e}_3 = (0,0,1)$. The fluid is bounded below by a rigid flat bed at $\{x_3 = -1\}$, while the upper boundary—described by the graph of $\eta$—represents a free interface with the air (modeled as a vacuum). A key feature of this free-boundary problem is that the free surface elevation $\eta$ is an unknown of the system. For solitary waves, we require $\eta(t, x_1, x_2) \to 0$ as $r = \sqrt{x_1^2 + x_2^2} \to \infty$, so the asymptotic fluid depth is normalized to 1.

	The boundary conditions for the system are specified as follows. At the rigid bottom \(x_3 = -1\), the kinematic condition of impermeability requires
	\begin{equation}
		\label{eq:impermeability}
		\mathbf{u} \cdot \mathbf{n} = 0,
	\end{equation}
	where \(\mathbf{n}\) is the unit outward normal to the fluid domain. On the free surface \(x_3 = \eta(t, x_1, x_2)\), the kinematic boundary condition is given by
	\begin{equation}
		\label{eq:kinematic}
		\partial_t \eta = \sqrt{1 + (\partial_{x_1} \eta)^2 + (\partial_{x_2} \eta)^2}  \, \mathbf{u} \cdot \mathbf{n},
	\end{equation}
	which states that the surface moves with the fluid. Dynamically, the pressure at the interface satisfies the Young–Laplace law:
	\begin{equation}
		\label{eq:pressure}
		P = \sigma \, \mathrm{div} \left( \frac{\nabla_{\!\perp} \eta}{\sqrt{1 + |\nabla_{\!\perp} \eta|^2}} \right),
	\end{equation}
	where \(\nabla_{\!\perp} = (\partial_{x_1}, \partial_{x_2})\) and \(\sigma > 0\) is the surface tension coefficient. Together with the gravitational acceleration \(g > 0\) in the bulk, the presence of surface tension at the interface leads to the terminology {\em capillary-gravity waves} for solutions of system \eqref{eq:euler1}--\eqref{eq:pressure}.

	Under the assumption of irrotational flow, the velocity field can be expressed as \(\mathbf{u} = \nabla \Phi\), and the system can be reformulated in terms of the velocity potential \(\Phi\). The resulting boundary value problem for the capillary-gravity water-wave system is:
	\begin{equation}
		\label{eq:water-wave-full}
		\begin{cases}
			\Delta \Phi = 0 & \text{in } \Omega(t), \\
			\partial_n \Phi = 0 & \text{on } x_3 = -1, \\
			\partial_t \eta = \partial_{x_3} \Phi - \nabla_{\perp} \eta \cdot \nabla_{\perp} \Phi & \text{on } x_3 = \eta(x_1, x_2, t), \\
			\partial_t \Phi = -\dfrac{1}{2} |\nabla \Phi|^2 - g \eta + \sigma \, \mathrm{div} \left( \dfrac{\nabla_{\perp} \eta}{\sqrt{1 + |\nabla_{\perp} \eta|^2}} \right) & \text{on } x_3 = \eta(x_1, x_2, t).
		\end{cases}
	\end{equation}
	
	A fundamental tool in the modern analysis of water waves is the Dirichlet-to-Neumann operator (DNO), which allows for a reduction of the problem's dimension by exploiting the harmonicity of the velocity potential. Defining the surface velocity potential as $\xi(x_1,x_2,t) = \Phi(x_1,x_2,\eta(x_1,x_2,t),t)$, the DNO $G(\eta)$ is defined by
	\begin{equation}
		G(\eta)\xi = \partial_{x_3}\Phi - \nabla_{\perp}\eta \cdot \nabla_{\perp}\Phi,
	\end{equation}
	where all quantities are evaluated at the free surface and $\nabla_\perp = (\partial_{x_1}, \partial_{x_2})$. In terms of the DNO, the water wave system \eqref{eq:water-wave-full} transforms into the Hamiltonian formulation \cite{Craig1,Craig2,Craig3}:
	\begin{equation}
		\label{eq:water-wave-bd}
		\begin{cases}
			\eta_t = G(\eta)\xi, \\
			\xi_t = \frac{1}{2(1+|\nabla_{\perp}\eta|^2)}\left[(G(\eta)\xi+\nabla_{\perp}\eta\cdot\nabla_{\perp}\xi)^2-|\nabla_{\perp}\xi|^2|\nabla_\perp\eta|^2-|\nabla_\perp\xi|^2\right] \\
			\quad\quad - g\eta + \sigma\operatorname{div} \left( \frac{\nabla_{\perp} \eta}{\sqrt{1 + |\nabla_{\perp} \eta|^2}} \right).
		\end{cases}
	\end{equation}
	This formulation reveals the Hamiltonian structure of the problem, with the energy functional given by
	\begin{equation}
		\label{1.energy-2}
		\mathcal{H}(\eta,\xi) = \int_{\mathbb{R}^2} \left( \frac12 \xi G(\eta)\xi + \frac12 g \eta^2 + \sigma \left( \sqrt{1+|\nabla_\perp \eta|^2} - 1 \right) \right) dx_1dx_2.
	\end{equation}
	The first term represents the kinetic energy, the second the gravitational potential energy, and the third the capillary energy due to surface tension.
	
	Due to the translation invariance of the system in the horizontal directions, there exist conserved momentum functionals. For waves propagating primarily in the $x_1$-direction, the relevant conserved quantity is the $x_1$-component of the momentum, given by
	\begin{equation}
		\label{eq:momentum}
		\mathcal{P}(\eta,\xi) = \int_{\mathbb{R}^2} \eta_{x_1} \xi  dx_1dx_2.
	\end{equation}
	This momentum functional plays a crucial role in both the existence theory and stability analysis of solitary waves \cite{Benjamin1974}.
	The DNO framework, along with the energy and momentum functionals, provides the mathematical foundation for both the existence theory and our subsequent stability analysis.

	A complete global well-posedness theory for the 3D gravity–capillary wave system in finite depth remains an open problem. (Existing global results, such as those by Germain–Masmoudi–Shatah \cite{GMS2012} and Wu \cite{Wu2009, Wu2011}, assume infinite depth and zero surface tension.) In the finite-depth setting, local well-posedness in 3D has been established by Alazard–Burq–Zuily \cite{ABZ2011} in the presence of bottom topography or for infinite depth, and by Ming–Zhang–Zhang \cite{MZZ2012} over long time scales.

	The existence theory for three-dimensional traveling waves features two main classes: doubly periodic waves and fully localized solitary waves. Doubly periodic waves, also known as short-crested waves, were first rigorously established by Reeder and Shinbrot \cite{Reed} via bifurcation theory. Earlier formal computations were given by Fuchs \cite{f1952} and Sretenskiĭ \cite{s1953}. Subsequent developments include the variational approach of Craig and Nicholls \cite{Cra2000} and spatial dynamics methods by Groves and Mielke \cite{gm2001}, Groves and Haragus \cite{gh2003}, and Nilsson \cite{n2019}. Recent advances for doubly periodic waves with vorticity include the work of Lokharu, Seth, and Wahlén \cite{lsw2020}, Groves et al. \cite{gnpw2024}, and Seth, Varholm, and Wahlén \cite{svw2024}.

	For fully localized solitary waves, where $\eta(x_1, x_2) \to 0$ as $|(x_1, x_2)| \to \infty$, two principal approaches have been successful. Variational methods seek critical points of the energy functional \eqref{1.energy-2} subject to momentum constraints, as in the work of Groves and Sun \cite{Groves} using mountain-pass arguments for the augmented energy $E_c(\eta, \xi) = \mathcal{H}(\eta, \xi) - c\mathcal{P}(\eta, \xi)$, and Buffoni et al. \cite{BGSW2013, Buffoni2018} via constrained minimization. In contrast, the waves we study were constructed by Gui, et al. \cite{GLLWY} using a non-variational Lyapunov–Schmidt reduction method, yielding small-amplitude solitary waves with speed $c = 1/\sqrt{1+\varepsilon^2}$ for $\varepsilon > 0$ sufficiently small. For a comprehensive survey of recent progress in water waves, we refer to \cite{H}.

	The primary objective of this paper is to establish the orbital stability of these fully localized three-dimensional solitary waves from \cite{GLLWY}. Crucially, these waves are not obtained as energy minimizers under momentum constraints, but rather as critical points through a reduction procedure. Consequently, their stability cannot be deduced from standard variational principles \cite{Caz, Benjamin1974} that apply to constrained minimizers \cite{Buffoni2004, Buffoni2005, Buffoni2009, Groves2010, Groves2011, Groves2015}.

	Our approach adapts the energy–momentum framework pioneered by Benjamin \cite{Benjamin} and rigorously developed by Grillakis, Shatah, and Strauss (GSS) \cite{GSS1,GSS2}. This method constructs a Lyapunov functional from a carefully chosen combination of the energy $\mathcal{H}$ and momentum $\mathcal{P}$, with stability following from the convexity properties of this functional evaluated at the wave profile. A key precedent is Mielke's work \cite{Mielke} on conditional energetic stability, which adapted the GSS method to accommodate the mismatch between well-posedness and energy spaces. While our strategy is philosophically similar, a fundamental distinction lies in the object of study: the waves we consider are not variational minimizers but arise from a non-variational construction. Consequently, we develop a framework to directly address the orbital stability of these specific wave profiles, synthesizing the GSS method with a detailed analysis of their linearized dynamics.

	Our main result establishes the conditional orbital stability of the wave profiles $\bar{u}(c)=(\bar{\eta}(c), \bar{\xi}(c))$ from \cite{GLLWY}:

	\begin{theorem}
		Every solitary capillary-gravity water wave $\bar{u}(c)=(\bar{\eta}(c), \bar{\xi}(c))$ with wave speed $c = \frac{1}{\sqrt{1+\varepsilon^2}}$ for $\varepsilon \in (0,\varepsilon_0)$ is conditionally orbitally stable in the following sense: For every $R > 0$ and $\rho > 0$, there exists $\rho_0 > 0$ such that if 
		$u=(\eta,\xi) \colon [0,T) \to \mathscr{F}_R$ is a continuous solution of the Hamiltonian system \eqref{eq:water-wave-bd} that preserves the energy $\mathcal{H}$ and momentum $\mathcal{P}$, and if the initial data satisfy
		\[
		\|\eta(0) - \bar{\eta}(c)\|_{H^1(\mathbb{R}^2)} + \|\xi(0) - \bar{\xi}(c)\|_{H^{1/2}_{*}(\mathbb{R}^2)} < \rho_0,
		\]
		then for all $t \in [0,T)$,
		\begin{equation*}
			\inf_{(x_0,y_0) \in \mathbb{R}^2} \big( \|\eta(t, \cdot - (x_0,y_0)) - \bar{\eta}(c)\|_{H^1} + \|\xi(t, \cdot - (x_0,y_0)) - \bar{\xi}(c)\|_{H^{1/2}_{*}} \big) < \rho.
		\end{equation*}
	\end{theorem}
	This theorem affirms the stability of specific wave profiles from \cite{GLLWY}, as opposed to the stability of a set of minimizers. Our proof carefully accounts for the distinctive structure of these reduction-based waves while handling the technical challenges posed by the functional setting. The conservation of both energy $\mathcal{H}$ and momentum $\mathcal{P}$ plays an essential role in the analysis, allowing us to construct an appropriate Lyapunov functional despite the non-variational origin of the waves.

	The paper is organized as follows. Section 2 introduces the Hamiltonian formulation via the Dirichlet-to-Neumann operator and establishes the appropriate functional setting for our analysis. Section 3 reviews the construction and properties of the solitary waves from \cite{GLLWY}, with particular emphasis on their asymptotic behavior and the spectral structure of the linearized operator. Section 4 completes the proof of our main stability theorem by implementing a refined version of the GSS method within Mielke's framework, carefully addressing the challenges posed by the non-variational nature of our waves. Technical estimates and auxiliary results are collected in the Appendix.

	\section{Hamiltonian Formulation and Stability Foundations}
	
	This section presents the Hamiltonian formulation of the three-dimensional capillary-gravity water wave problem and derives the key variational identities that form the foundation of our stability analysis. We begin by recalling the fundamental structure of the water wave system in the fixed reference domain $\Omega=\{(x_1,x_2,x_3)\in \mathbb R^2\times (-1,0)\}$.

	Traveling wave solutions propagating with constant speed $c$ in the $x_1$-direction satisfy the steady system:
	\begin{equation}\label{1-1}
		\begin{cases}
			\begin{aligned}
				-c\partial_{1}\eta &= G(\eta) \xi, \\
				-c\partial_{1}\xi &= \frac{1}{2(1+|\nabla \eta|^2)} \left[ (G(\eta)\xi + \nabla\eta \cdot \nabla\xi)^2 - |\nabla\xi|^2 |\nabla\eta|^2 - |\nabla\xi|^2 \right] \\
				&\quad - g\eta + \sigma \operatorname{div} \left( \frac{\nabla \eta}{\sqrt{1 + |\nabla \eta|^2}} \right),
			\end{aligned}
		\end{cases}
	\end{equation}
	where $\eta, \xi: \mathbb{R}^2 \to \mathbb{R}$ denote the free surface and velocity potential, respectively, and $\nabla = (\partial_{x_1}, \partial_{x_2})$.

	A fundamental observation is that the pair $(\eta,\xi)$ constitutes a solution of \eqref{1-1} if and only if it is a critical point of the augmented Hamiltonian functional $\mathcal{H}_c:=\mathcal{H}+c\mathcal{P}$. This variational characterization provides the essential link between the existence theory and our stability analysis. The energy and momentum functionals are defined by
	\[
	\mathcal H(u)=\mathcal{K}(u)+\mathcal{V}(u),\quad  \mathcal{P}(u)=\int_{\mathbb{R}^2}\partial_{x_1}\eta \xi dx_1 dx_2,
	\]
	with kinetic and potential energy components given respectively by
	\[
	\mathcal{K}(u)=\int_{\mathbb{R}^2} \frac{1}{2} \xi(G(\eta) \xi)dx_1 dx_2,\quad 
	\mathcal{V}(u)=\int_{\mathbb{R}^2}\left[\frac{1}{2} g \eta^2+\sigma\left(\sqrt{1+\left|\nabla \eta\right|^2}-1\right)\right] d x_1 dx_2.
	\]
	The time-dependent water wave system admits the canonical Hamiltonian structure
	\begin{equation*}
		\partial_t u=\nabla_{u}^{\perp}\mathcal{H}(u),\quad u=(\eta,\xi),\quad \nabla_{u}^{\perp}=(\partial_{\xi},-\partial_{\eta}),
	\end{equation*}
	which reveals the underlying symplectic geometry of the fluid dynamics.

	For the stability analysis of traveling waves, it proves advantageous to work with the reduced augmented potential obtained by minimizing the Hamiltonian with respect to the surface potential variable. This reduction transforms the search for critical points into a minimization problem, significantly simplifying the subsequent analysis. Specifically, we define
	\begin{align*}
		\mathcal{V}^{aug}_{c}(\eta):&=\min\{ \mathcal{H}_c(\eta,\xi) \mid \xi \in H^{1/2}_{*}(\mathbb R^2)\}\\
		&=\mathcal{V}(\eta)-\frac{c^2}{2}\int_{\mathbb{R}^2}\partial_{x_1}\eta [G(\eta)^{-1}\partial_{x_1}\eta]dx_1dx_2.
	\end{align*}
	The first and second variations of the potential energy component $\mathcal{V}$ play a crucial role in our stability criteria and are computed directly as
	\begin{equation*}
		D\mathcal{V}(\eta)[v]=\int_{\mathbb R^2}\left[g\eta v+\sigma\frac{\nabla\eta \cdot \nabla v}{\sqrt{1+|\nabla \eta|^2}} \right]dx_1dx_2
	\end{equation*}
	and
	\begin{equation*}
		D^2\mathcal{V}(\eta)[v,v]=\int_{\mathbb R^2} \left[ g v^2 + \sigma\Big( \frac{|\nabla v|^2}{\sqrt{1+|\nabla \eta|^2}} - \frac{(\nabla \eta \cdot \nabla v)^2}{(1+|\nabla \eta|^2)^{3/2}} \Big) \right] dx_1dx_2.
	\end{equation*}
	We now derive the second variation of the augmented potential, which governs the linear stability properties of traveling wave solutions. This computation forms the analytical core of our approach.
	\begin{lemma}\label{2-3}
		If $u = (\eta, \xi)$ is a critical point of $\mathcal{H}_c$, then for all test functions $v$ the second variation of the augmented potential satisfies the identity:
		\begin{equation*}
			D^2 \mathcal{V}_c^{\text{aug}}(\eta)[v, v] = D^2 \mathcal{V}(\eta)[v, v] + \frac{1}{2} \left\langle D^2 G(\eta)[v, v] \xi, \xi \right\rangle - \left\langle \mathcal{L}_u v, G(\eta)^{-1} \mathcal{L}_u v \right\rangle,
		\end{equation*}
		where the linearized transport operator is given by
		\[
		\mathcal{L}_u v = c \partial_{x_1} v + D G(\eta)[v] \xi.
		\]
	\end{lemma}
	
	\begin{proof}
		We begin with the expression for the second variation of the augmented Hamiltonian:
		\begin{equation}\label{2.3}
			D^2\mathcal{H}_c(u)\left[\binom{v}{w}, \binom{v}{w}\right] = D^2 \mathcal{V}(\eta)[v, v] + \frac{1}{2} \left\langle D^2 G(\eta)[v, v] \xi, \xi \right\rangle + \langle G(\eta) w, w \rangle + 2 \langle \mathcal{L}_u v, w \rangle.
		\end{equation}
		To obtain the second variation of the augmented potential, we minimize this expression with respect to $w$. Define the quadratic form in $w$ as
		\[
		Q(w) = \langle G(\eta) w, w \rangle + 2 \langle \mathcal{L}_u v, w \rangle.
		\]
		The Dirichlet-to-Neumann operator $G(\eta)$ is uniformly positive definite on $H^{1/2}_*(\mathbb{R}^2)$ (see Lemma \ref{DNO bound}), ensuring that $Q(w)$ is strictly convex and coercive. The unique minimizer is found by solving the Euler-Lagrange equation $DQ(w) = 0$, which yields
		\[
		G(\eta) w = -\mathcal{L}_u v.
		\]
		Substituting this optimal $w$ back into equation \eqref{2.3} and observing that the cross term becomes
		\[
		2\langle \mathcal{L}_u v, w \rangle = -2\langle \mathcal{L}_u v, G(\eta)^{-1}\mathcal{L}_u v \rangle,
		\]
		while the quadratic term becomes
		\[
		\langle G(\eta) w, w \rangle = \langle \mathcal{L}_u v, G(\eta)^{-1}\mathcal{L}_u v \rangle,
		\]
		we obtain the desired identity and complete the proof.
	\end{proof}

	To simplify the expression for $D^2\mathcal{V}^{\text{aug}}_c$, we aim to express it in terms of pointwise functions as much as possible. This pointwise representation will facilitate our subsequent spectral analysis and stability arguments. We begin by computing explicit expressions for the first and second variations of the Dirichlet-to-Neumann operator, which play a central role in the second variation formula.

	\begin{lemma}\label{DG}
		The first variation of the Dirichlet-to-Neumann operator satisfies
		\begin{equation}\label{1DG}
			\begin{aligned}
				\langle \psi, {DG}(\eta)[v] \phi\rangle 
				& = \left\langle v, a_0(\eta, \phi, \psi)\right\rangle \\
				& = \left\langle v, a_1(\eta, \phi)\cdot \nabla \psi + a_2(\eta, \phi) G(\eta) \psi\right\rangle,
			\end{aligned}
		\end{equation}
		where the coefficient functions are given by
		\begin{align*}
			a_0(\eta,\phi,\psi) & = \nabla \phi\cdot \nabla \psi - \frac{1}{1+|\nabla \eta|^2}\left[ (G(\eta)\phi + \nabla \eta\cdot \nabla \phi)\cdot (G(\eta)\psi + \nabla \eta\cdot \nabla \psi) \right], \\
			a_1(\eta,\phi) & = \nabla \phi - \frac{[G(\eta)\phi + \nabla \eta \cdot \nabla \phi]}{1+|\nabla \eta|^2} \nabla \eta, \\
			a_2(\eta,\phi) & = -\frac{[G(\eta)\phi + \nabla \eta \cdot \nabla \phi]}{1+|\nabla \eta|^2}.
		\end{align*}
		For the second variation we have the representation
		\begin{equation*}
			\langle \phi, {D^2G}(\eta)[v,v] \phi \rangle = \langle v, a_3 v\rangle + 2\langle a_2v, G(\eta)(a_2 v)\rangle,
		\end{equation*}
		where $a_3 = -2a_2 \operatorname{div}(a_1)$.
	\end{lemma}

	\begin{proof}
		We derive the first derivative of the Dirichlet–Neumann operator $G(\eta)$ associated with the domain
		$$\Omega(\eta) = \{(x_1, x_2, x_3) \in \mathbb{R}^3 \mid (x_1, x_2) \in \mathbb{R}^2, \, -1 < x_3 < \eta(x_1, x_2)\},$$
		where $\inf \eta > -1$, and $\eta \in W^{1,\infty}(\mathbb{R}^2)$. The operator $G(\eta): H^1(\mathbb{R}^2) \to L^2(\mathbb{R}^2)$ is defined as follows: for $\phi \in H^1(\mathbb{R}^2)$, let $\Phi$ solve the boundary-value problem
		$$\Delta \Phi = 0 \quad \text{in } \Omega(\eta), \quad \partial_{x_3} \Phi = 0 \quad \text{at } x_3 = -1, \quad \Phi = \phi \quad \text{at } x_3 = \eta(x_1, x_2).$$
		Then,
		$$[G(\eta)\phi](x_1, x_2) = \partial_{x_3} \Phi(x_1, x_2, \eta(x_1, x_2)) - \nabla_{(x_1, x_2)}\eta \cdot \nabla_{(x_1, x_2)}  \Phi(x_1, x_2, \eta(x_1, x_2)).$$
		To compute the Fréchet derivative $DG(\eta)[v]$ in the direction $v \in W^{1,\infty}(\mathbb{R}^2)$ with $\inf (\eta + \varepsilon v) > -1$ for small $\varepsilon > 0$, we use the bilinear form representation. For test functions $\phi, \psi \in H^1(\mathbb{R}^2)$, let $\Phi$ and $\Psi$ be the respective harmonic extensions (lifts) as above.  The operator satisfies
		$$\langle \psi, G(\eta)\phi \rangle = \int_{\Omega(\eta)} \sum_{i=1}^3 \partial_{x_i} \Phi \cdot \partial_{x_i} \Psi \, dx_1dx_2dx_3,$$
		where $\langle \cdot, \cdot \rangle = \int_{\mathbb{R}^2} \cdot \, dx_1dx_2$ denotes the $L^2(\mathbb{R}^2)$ inner product. This follows from Green's identity, noting the Neumann condition at the bottom $x_3 = -1$ and the relation $G(\eta)\phi = \sqrt{1 + |\nabla \eta|^2} \, \partial_n \Phi$, where $\partial_n$ is the outward normal derivative. Consider the perturbed surface $\eta_\varepsilon = \eta + \varepsilon v$ and domain $\Omega_\varepsilon = \Omega(\eta_\varepsilon)$. Let $\Phi_\varepsilon$ and $\Psi_\varepsilon$ be the lifts of $\phi$ and $\psi$ in $\Omega_\varepsilon$, respectively. Then,
		$$\langle \psi, G(\eta_\varepsilon)\phi \rangle = \int_{\Omega_\varepsilon} \sum_{i=1}^3 \partial_{x_i}  \Phi_\varepsilon \cdot \partial_{x_i} \Psi_\varepsilon \, dx_1dx_2dx_3.$$
		Differentiating with respect to $\varepsilon$ at $\varepsilon = 0$ gives 
		$$\langle \psi, DG(\eta)[v]\phi \rangle = \frac{\mathrm{d}}{\mathrm{d}\varepsilon} \bigg|_{\varepsilon=0} \int_{\Omega_\varepsilon} \sum_{i=1}^3\partial_{x_i} \Phi_\varepsilon \cdot \partial_{x_i} \Psi_\varepsilon \, dx_1dx_2dx_3.$$
		By the Reynolds transport theorem for the varying domain (with the bottom fixed and only the top boundary perturbed vertically by $\varepsilon v$),
		$$\frac{\mathrm{d}}{\mathrm{d}\varepsilon} \int_{\Omega_\varepsilon} F_\varepsilon \, dx_1dx_2dx_3 = \int_{\Omega} \partial_\varepsilon F_\varepsilon \, dx_1dx_2dx_3 + \int_{\partial \Omega^\text{top}} F_0 (V \cdot n) \, {d}S,$$
		where $F_\varepsilon = \sum_{i=1}^3\partial_{x_i} \Phi_\varepsilon \cdot \partial_{x_i} \Psi_\varepsilon$, $\Omega = \Omega(\eta)$, $V = (0,0, v)$ is the boundary velocity field at the top, $n$ is the outward unit normal at the top boundary $x_3 = \eta(x_1,x_2)$, i.e., $n = (-\nabla\eta, 1)/\sqrt{1 + |\nabla \eta|^2}$, and ${d}S = \sqrt{1 + |\nabla \eta|^2} \, dx_1dx_2$.
		The perturbation yields $V \cdot n = v / \sqrt{1 + |\nabla \eta|^2}$, so
		$$\int_{\partial \Omega^\text{top}} F_0 (V \cdot n) \, dS = \int_{\mathbb{R}^2} [\sum_{i=1}^3 \partial_{x_i} \Phi \cdot \partial_{x_i} \Psi]_{x_3=\eta(x_1,x_2)} v(x_1,x_2) \, dx_1dx_2.$$
		The field derivative term is
		$$\partial_\varepsilon F_\varepsilon \big|_{\varepsilon=0} = \sum_{i=1}^3 \partial_{x_i} \dot{\Phi} \cdot \partial_{x_i} \Psi + \partial_{x_i} \Phi \cdot \partial_{x_i} \dot{\Psi},$$
		where $\dot{\Phi} = \partial_\varepsilon \Phi_\varepsilon |_{\varepsilon=0}$ and $\dot{\Psi} = \partial_\varepsilon \Psi_\varepsilon |_{\varepsilon=0}$ are the shape derivatives. 
		Since $\Phi_\varepsilon(x_1,x_2, \eta_\varepsilon(x_1,x_2)) = \phi(x_1,x_2)$  for all  $(x_1,x_2) \in \mathbb{R}^2$, it follows that  
		the shape derivatives satisfy $\Delta \dot{\Phi} = 0$ in $\Omega$, $\partial_{x_3} \dot{\Phi} = 0$ at $x_3 = -1$, and at $x_3 = \eta(x_1, x_2)$,
		$$\dot{\Phi}(x_1, x_2, \eta(x_1, x_2)) = -v(x_1, x_2) \partial_{x_3} \Phi(x_1, x_2, \eta(x_1, x_2)),$$
		with an analogous boundary condition for $\dot{\Psi}$:
		$$\dot{\Psi}(x_1, x_2, \eta(x_1, x_2)) = -v(x_1, x_2) \partial_{x_3} \Psi(x_1, x_2, \eta(x_1, x_2)).$$
		Using the bilinear form identity,
		\begin{align*}
			\int_{\Omega} \sum_{i=1}^3\partial_{x_i} \dot{\Phi} \cdot \partial_{x_i} \Psi \, dx_1dx_2dx_3 =
			\langle \dot{\Phi}|_{x_3 =\eta}, G(\eta)\psi \rangle 
			= -\int_{\mathbb{R}^2} v (\partial_{x_3} \Phi) [G(\eta)\psi] \, dx_1dx_2,
		\end{align*}
		and similarly,
		\begin{align*}
			\int_{\Omega} \sum_{i=1}^3\partial_{x_i} \Phi \cdot \partial_{x_i} \dot{\Psi} \, dx_1dx_2dx_3 =\langle \phi, G(\eta) [-v \partial_{x_3} \Psi] \rangle = -\int_{\mathbb{R}^2} v (\partial_{x_3} \Psi) [G(\eta)\phi] \, dx_1dx_2,
		\end{align*}
		where the last step uses self-adjointness of $G(\eta)$.
		Thus,
		\begin{align*}
			\langle \psi, DG(\eta)[v]\phi \rangle &= -\int_{\mathbb{R}^2} v (\partial_{x_3} \Phi) [G(\eta)\psi] \, dx_1dx_2 - \int_{\mathbb{R}^2} v (\partial_{x_3} \Psi) [G(\eta)\phi] \, dx_1dx_2\\
			&\quad + \int_{\mathbb{R}^2} v (\nabla_{(x_1, x_2)} \Phi \nabla_{(x_1, x_2)} \Psi + \partial_{x_3} \Phi \partial_{x_3} \Psi) \, dx_1dx_2,
		\end{align*}
		with all boundary terms evaluated at $x_3 = \eta(x_1, x_2)$.
		At the boundary, $\nabla_{(x_1, x_2)} \Phi = \nabla\phi - \nabla \eta \partial_{x_3} \Phi$ and $\nabla_{(x_1, x_2)} \Psi = \nabla\psi - \nabla \eta \partial_{x_3} \Psi$, so
		\begin{align*}
			&\nabla_{(x_1, x_2)} \Phi \nabla_{(x_1, x_2)} \Psi + \partial_{x_3} \Phi \partial_{x_3} \Psi\\
			&= (\nabla \phi - \nabla \eta  \partial_{x_3} \Phi)(\nabla\psi - \nabla \eta \partial_{x_3} \Psi) + \partial_{x_3} \Phi \partial_{x_3} \Psi\\
			&= \nabla\phi \cdot \nabla \psi - \nabla\eta (\partial_{x_3} \Phi \nabla \psi + \partial_{x_3} \Psi \nabla\phi) + (1 + |\nabla \eta| ^2) \partial_{x_3} \Phi \partial_{x_3} \Psi.
		\end{align*}
		Additionally, $\partial_{x_3} \Phi = [G(\eta)\phi + \nabla \eta \cdot \nabla\phi] / [1 + |\nabla \eta|^2]$ and $\partial_{x_3} \Psi = [G(\eta)\psi + \nabla \eta \cdot \nabla \psi] / [1 + |\nabla \eta|^2]$. Substituting these and collecting terms yields
		\begin{align*}
			\langle \psi, DG(\eta)[v]\phi \rangle &= \int_{\mathbb{R}^2} v \nabla \phi \cdot \nabla \psi dx_1dx_2\\
			&\quad-\int_{\mathbb{R}^2}\frac{v}{1+|\nabla \eta|^2}(G(\eta)\phi+\nabla \eta\cdot \nabla \phi) (G(\eta)\psi+\nabla \eta\cdot \nabla \psi) dx_1dx_2.
		\end{align*}
		Hence,
		$$\langle \psi, DG(\eta)[v]\phi \rangle = \langle v, a_0(\eta, \phi, \psi) \rangle,$$
		where
		$$a_0(\eta, \phi, \psi) = \nabla \phi \cdot \nabla \psi- \frac{1}{1+|\nabla \eta|^2} (G(\eta)\phi+\nabla \eta\cdot \nabla \phi) (G(\eta)\psi+\nabla \eta\cdot \nabla \psi).$$
		This can be rewritten as
		$$a_0(\eta, \phi, \psi) = a_1(\eta, \phi) \nabla \psi + a_2(\eta, \phi) [G(\eta)\psi],$$
		with
		\begin{align*}
			a_1(\eta,\phi)=\nabla \phi-\frac{[G(\eta)\phi +\nabla \eta \cdot \nabla \phi]}{1+|\nabla \eta|^2}\nabla \eta, \quad
			a_2(\eta,\phi)=-\frac{[G(\eta)\phi +\nabla \eta \cdot \nabla \phi]}{1+|\nabla \eta|^2}.
		\end{align*}
		To compute the second derivative, differentiate the first derivative with respect to $\eta$ along the direction $v$.
		Given that $a_1(\eta,\phi) = \nabla \phi + a_2(\eta,\phi) \nabla \eta$, we obtain
		\begin{align*}
			\langle \phi, D^2G(\eta)[v,v]\phi \rangle
			&= \langle v, a_2 \nabla v \cdot \nabla \phi + Da_2[v] \nabla \eta \cdot \nabla \phi + Da_2[v] G(\eta) \phi + a_2 DG(\eta)[v] \phi \rangle,
		\end{align*}
		where
		\begin{align*}
			Da_2[v] = -\frac{2\nabla \eta \cdot \nabla v}{1 + |\nabla \eta|^2}a_2 - \frac{DG(\eta)[v]\phi + \nabla v \cdot \nabla \phi}{1 + |\nabla \eta|^2}.
		\end{align*}
		It follows that
		\begin{align*}
			\langle \phi, D^2G(\eta)[v,v]\phi \rangle
			&= \langle v, a_2 \nabla v \cdot \nabla \phi - Da_2[v]\left(a_2(1 + |\nabla \eta|^2)\right) + a_2 DG(\eta)[v] \phi \rangle \\
			&= \langle v, 2a_2 \nabla v \cdot \nabla \phi + 2a_2^2 \nabla \eta \cdot \nabla v + 2a_2 DG(\eta)[v] \phi \rangle \\
			&= \langle v, 2a_2 a_1 \cdot \nabla v + 2a_2 DG(\eta)[v] \phi \rangle,
		\end{align*}
		with
		\begin{align*}
			\langle v, 2a_2 DG(\eta)[v] \phi \rangle
			&= 2 \langle a_2 v, DG(\eta)[v] \phi \rangle \\
			&= 2 \langle v, a_1 \cdot \nabla(a_2 v) + a_2 G(\eta)(a_2 v) \rangle \\
			&= 2 \langle v, a_1 \cdot (\nabla a_2) v + a_2 a_1 \cdot \nabla v \rangle + 2 \langle v, a_2 G(\eta)(a_2 v) \rangle.
		\end{align*}
		Therefore,
		\begin{align*}
			\langle \phi, D^2G(\eta)[v,v]\phi \rangle
			&= \langle v, 4a_2 a_1 \cdot \nabla v + 2 (a_1 \cdot \nabla a_2) v \rangle + 2 \langle v, a_2 G(\eta)(a_2 v) \rangle \\
			&= \langle v, -2 \nabla \cdot (a_2 a_1) v + 2 (a_1 \cdot \nabla a_2) v \rangle + 2 \langle v, a_2 G(\eta)(a_2 v) \rangle \\
			&= \langle v, -2a_2 (\nabla \cdot a_1) v \rangle + 2 \langle v, a_2 G(\eta)(a_2 v) \rangle.
		\end{align*}
		This completes the proof.
	\end{proof}

	With the explicit expressions for the variations of the Dirichlet-to-Neumann operator established in Lemma \ref{DG}, we are now prepared to derive the pointwise representation of the second variation of the augmented potential. This representation will reveal the underlying structure of the linearized operator and provide the foundation for our spectral analysis. The following theorem gives the precise form of $D^2\mathcal{V}^{\mathrm{aug}}_c(\eta)$ in terms of local coefficients and a nonlocal, self-adjoint operator.

	\begin{theorem}\label{DV}
		Let \( u = (\eta, \xi) \) be a critical point of \( \mathcal{H}_c \). Then the second variation of the augmented potential \( \mathcal{V}^{\mathrm{aug}}_c \) at \( \eta \) takes the form
		\begin{align*}
			D^2\mathcal{V}^{\mathrm{aug}}_c(\eta)[v, v] 
			= \int_{\mathbb{R}^2} \left[ - v \, \mathrm{div}( \alpha \nabla v) + \beta v^2 - v (M_{\eta} v) \right]  dx_1dx_2,
		\end{align*}
		where the coefficient tensor $\alpha = (\alpha_{ij})_{2\times 2}$ has components
		\[
		\alpha_{ij} = \frac{\sigma}{\left(1 + |\nabla \eta|^2\right)^{3/2}} 
		\left[ 
		\left(1 + |\nabla \eta|^2\right) \delta_{ij} - \partial_i \eta \, \partial_j \eta 
		\right],
		\quad i,j = 1,2.
		\]
		Here, $\delta_{ij}$ denotes the Kronecker delta, and $\partial_i \eta$ represents the partial derivative of $\eta$ with respect to the $i$-th spatial coordinate.
		The scalar coefficient \( \beta \) is given by
		\[
		\beta = g - \rho_1 \cdot \nabla \rho_2 + c \partial_{x_1} \rho_2,
		\]
		where
		\[
		\rho_1 = \nabla \xi - \frac{\nabla \eta \cdot \nabla \xi - c \partial_{x_1} \eta}{1 + |\nabla \eta|^2} \nabla \eta, 
		\quad 
		\rho_2 = \frac{c \partial_{x_1} \eta - \nabla \eta \cdot \nabla \xi}{1 + |\nabla \eta|^2}.
		\]
		Furthermore, \( M_{\eta} \) denotes a self-adjoint positive operator on \( L^2(\mathbb{R}^2) \) defined by
		\[
		M_{\eta} v = - \left( c \partial_{x_1} - \rho_1 \cdot \nabla \right) 
		\left[ G(\eta)^{-1} \left( c \partial_{x_1} v - \mathrm{div}(v \rho_1) \right) \right].
		\]
	\end{theorem}

	\begin{proof}
		By Lemma \ref{2-3}, we have
		\begin{equation*}
			D^2\mathcal{V}_c^{\text{aug}}(\eta)[v,v] = D^2\mathcal{V}(\eta)[v,v] + \frac{1}{2}\langle D^2G(\eta)[v,v]\xi,\xi\rangle - \langle\mathcal{L}_u v, G(\eta)^{-1}\mathcal{L}_u v\rangle.
		\end{equation*}
		We compute each term separately. 
		First, the second variation of the potential energy gives
		\begin{equation*}
			D^2\mathcal{V}(\eta)[v,v]=\int_{\mathbb R^2} \left[ g v^2 + \sigma\left( \frac{|\nabla v|^2}{\sqrt{1+|\nabla \eta|^2}} - \frac{(\nabla \eta \cdot \nabla v)^2}{(1+|\nabla \eta|^2)^{3/2}} \right) \right] dx_1dx_2.
		\end{equation*}
		Since 
		\begin{align*}
			\sigma\left( \frac{|\nabla v|^2}{\sqrt{1+|\nabla \eta|^2}} - \frac{(\nabla \eta \cdot \nabla v)^2}{(1+|\nabla \eta|^2)^{3/2}} \right) &= \sigma \frac{ |\nabla v|^2 (1 + |\nabla \eta|^2) - (\nabla \eta \cdot \nabla v)^2 }{(1+|\nabla \eta|^2)^{3/2}} \\
			&= \nabla v \cdot (\alpha \nabla v),
		\end{align*}
		where $\alpha = (\alpha_{ij})_{2\times 2}$ is defined by
		\[
		\alpha_{ij} = \frac{\sigma}{\left(1 + |\nabla \eta|^2\right)^{3/2}} 
		\left[ 
		\left(1 + |\nabla \eta|^2\right) \delta_{ij} - \partial_i \eta \, \partial_j \eta 
		\right],
		\quad i,j = 1,2,
		\]
		with $\delta_{ij}$ denoting the Kronecker delta, we get
		\begin{equation*}
			D^2\mathcal{V}(\eta)[v,v]=\int_{\mathbb R^2}\left[ g v^2 - v \operatorname{div}( \alpha \nabla v) \right] dx_1dx_2.
		\end{equation*}
		For the second term, applying Lemma \ref{DG} with $w = v$ and $\phi = \xi$, and using the critical point condition $G(\eta)\xi = -c\partial_{x_1}\eta$, we obtain
		\begin{equation*}
			\frac{1}{2}\langle D^2G(\eta)[v,v]\xi,\xi\rangle = \frac{1}{2}\langle v, \rho_3 v\rangle + \langle \rho_2 v, G(\eta)(\rho_2 v)\rangle,
		\end{equation*}
		where $\rho_3 = -2\rho_2\text{div}(\rho_1)$ with $\rho_1=\nabla \xi-\frac{[\nabla \eta \cdot \nabla \xi- c\partial_{x_1}\eta]}{1+|\nabla \eta|^2} \nabla \eta$ and $\rho_2=\frac{c\partial_{x_1}\eta-\nabla \eta \cdot \nabla \xi}{1+|\nabla \eta|^2}$.
		For the third term, we first compute $\mathcal{L}_u v $. From equation \eqref{1DG} with $\hat{v}$ as a test function, we have
		\begin{equation*}
			\langle \hat{v}, DG(\eta)[v]\xi\rangle = \langle v, \rho_1\cdot\nabla\hat{v} + \rho_2 G(\eta)\hat{v}\rangle.
		\end{equation*}
		This implies the distributional identity
		\begin{equation*}
			DG(\eta)[v]\xi = -\text{div}(v\rho_1) + G(\eta)(\rho_2 v),
		\end{equation*}
		and therefore
		\begin{equation*}
			\mathcal{L}_u v = c\partial_{x_1}v + DG(\eta)[v]\xi = c\partial_{x_1}v - \text{div}(v\rho_1) + G(\eta)(\rho_2 v).
		\end{equation*}
		Now we compute the inner product
		\begin{align*}
			&\langle\mathcal{L}_u v, G(\eta)^{-1}\mathcal{L}_u v\rangle \\
			&= \langle c\partial_{x_1}v - \text{div}(v\rho_1) + G(\eta)(\rho_2 v), G(\eta)^{-1}(c\partial_{x_1}v - \text{div}(v\rho_1) + G(\eta)(\rho_2 v))\rangle.
		\end{align*}
		Let $U = c\partial_{x_1}v - \text{div}(v\rho_1)$ and $V = G(\eta)(\rho_2 v)$. Then
		\begin{align*}
			\langle U + V, G(\eta)^{-1}(U + V)\rangle &= \langle U, G(\eta)^{-1}U\rangle + 2\langle U, \rho_2 v\rangle + \langle V, G(\eta)^{-1}V\rangle \\
			&= \langle U, G(\eta)^{-1}U\rangle + 2\langle U, \rho_2 v\rangle + \langle G(\eta)(\rho_2 v), \rho_2 v\rangle.
		\end{align*}
		We now compute the cross term $2\langle U, \rho_2 v\rangle$:
		\begin{align*}
			2\langle U, \rho_2 v\rangle &= 2\langle c\partial_{x_1}v - \text{div}(v\rho_1), \rho_2 v\rangle \\
			&= 2\int_{\mathbb{R}^2} (c\partial_{x_1}v)(\rho_2 v) dx_1dx_2 - 2\int_{\mathbb{R}^2} \text{div}(v\rho_1)(\rho_2 v) dx_1dx_2.
		\end{align*}
		Using integration by parts,
		\begin{align*}
			\int_{\mathbb{R}^2} (c\partial_{x_1}v)(\rho_2 v) dx_1dx_2 &= -\frac{1}{2}\int_{\mathbb{R}^2} c(\partial_{x_1}\rho_2)v^2 dx_1dx_2, \\
			\int_{\mathbb{R}^2} \text{div}(v\rho_1)(\rho_2 v) dx_1dx_2 &= -\int_{\mathbb{R}^2} v\rho_1\cdot\nabla(\rho_2 v) dx_1dx_2 \\
			&= -\int_{\mathbb{R}^2} [\rho_1\cdot\nabla\rho_2]v^2 dx_1dx_2 - \int_{\mathbb{R}^2} [\rho_1\cdot\nabla v](\rho_2 v) dx_1dx_2.
		\end{align*}
		For the last term,
		\begin{align*}
			\int_{\mathbb{R}^2} [\rho_1\cdot\nabla v](\rho_2 v) dx_1dx_2 &= \frac{1}{2}\int_{\mathbb{R}^2} \rho_1\rho_2\cdot\nabla(v^2) dx_1dx_2 \\
			&= -\frac{1}{2}\int_{\mathbb{R}^2} \text{div}(\rho_1\rho_2)v^2 dx_1dx_2.
		\end{align*}
		Combining these results,
		\begin{align*}
			2\langle U, \rho_2 v\rangle &= \int_{\mathbb{R}^2} [-c\partial_{x_1}\rho_2 + 2\rho_1\cdot\nabla\rho_2 - \text{div}(\rho_1\rho_2)]v^2 dx_1dx_2.
		\end{align*}
		Now, substituting all terms into the original expression,
		\begin{align*}
			&D^2\mathcal{V}_c^{\text{aug}}(\eta)[v,v] \\
			&= \int_{\mathbb{R}^2} gv^2- v \operatorname{div}( \alpha \nabla v) dx_1dx_2 + \frac{1}{2}\langle v, \rho_3 v\rangle + \langle \rho_2 v, G(\eta)(\rho_2 v)\rangle \\
			&\quad - \langle U, G(\eta)^{-1}U\rangle - 2\langle U, \rho_2 v\rangle - \langle G(\eta)(\rho_2 v), \rho_2 v\rangle.
		\end{align*}
		Note that $\langle \rho_2 v, G(\eta)(\rho_2 v)\rangle$ and $-\langle G(\eta)(\rho_2 v), \rho_2 v\rangle$ cancel. Also, $$\frac{1}{2}\langle v, \rho_3 v\rangle = -\langle v, \rho_2\text{div}(\rho_1) v\rangle.$$ Therefore, we get
		\begin{align*}
			&D^2\mathcal{V}_c^{\text{aug}}(\eta)[v,v] \\
			&= \int_{\mathbb{R}^2} gv^2- v \operatorname{div}( \alpha \nabla v) dx_1dx_2 - \langle v, \rho_2\text{div}(\rho_1) v\rangle - \langle U, G(\eta)^{-1}U\rangle \\
			&\quad - \int_{\mathbb{R}^2} [-c\partial_{x_1}\rho_2 + 2\rho_1\cdot\nabla\rho_2 - \text{div}(\rho_1\rho_2)]v^2 dx_1dx_2.
		\end{align*}
		Using $\text{div}(\rho_1\rho_2) = \rho_2\text{div}(\rho_1) + \rho_1\cdot\nabla\rho_2$, we simplify:
		\begin{align*}
			D^2\mathcal{V}_c^{\text{aug}}(\eta)[v,v] &= \int_{\mathbb{R}^2} - v \operatorname{div}( \alpha \nabla v) + \beta v^2 dx_1dx_2 - \langle U, G(\eta)^{-1}U\rangle,
		\end{align*}
		where $\beta = g+c\partial_{x_1}\rho_2 - \rho_1\cdot\nabla\rho_2$.
		Finally, we observe that
		\begin{align*}
			\langle U, G(\eta)^{-1}U\rangle &= \langle c\partial_{x_1}v - \text{div}(v\rho_1), G(\eta)^{-1}(c\partial_{x_1}v - \text{div}(v\rho_1))\rangle \\
			&= \langle v, M_\eta v\rangle,
		\end{align*}
		where $M_\eta$ is the self-adjoint positive operator defined by:
		\begin{equation*}
			M_\eta v = -\left(c\partial_{x_1} - \rho_1\cdot\nabla\right)\left(G(\eta)^{-1}(c\partial_{x_1}v - \text{div}(v\rho_1))\right).
		\end{equation*}
		This completes the proof.
	\end{proof}

	The quadratic form representation obtained in Theorem \ref{DV} naturally gives rise to an associated self-adjoint operator that governs the linearized dynamics around the traveling wave solution. This operator will play a central role in our stability analysis, as its spectral properties determine the linear stability of the wave.
	\begin{equation}\label{A}
		\begin{cases}
			A(\eta)v = -\operatorname{div}(\alpha \nabla v) + \beta v - M_{\eta}v,\\
			\mathcal{D}(A(\eta)) = H^2(\mathbb R^2).
		\end{cases}
	\end{equation}
	
	The operator $A(\eta)$ combines local elliptic terms with the nonlocal operator $M_{\eta}$, reflecting the mixed local-nonlocal character of the water wave problem. Its domain is naturally chosen as $H^2(\mathbb{R}^2)$ to ensure the proper treatment of both the second-order elliptic part and the nonlocal component.
	
	We now establish the fundamental spectral properties of this operator, which are crucial for understanding the stability characteristics of solitary waves.
	\begin{theorem}\label{thm2.1}
		Assume that $\eta \in H^s(\mathbb{R}^2)$ with $s > \frac{5}{2}$ and that $\inf\{\eta(x) \mid x \in \mathbb{R}^2\} > -h$. Then the operator $A(\eta)$ defined in \eqref{A} is self-adjoint on $L^2(\mathbb{R}^2)$ with domain $\mathcal{D}(A(\eta)) := H^2(\mathbb{R}^2)$ and has the continuous spectrum $[\sigma_*, +\infty)$ with
		\[
		\sigma_* = 
		\begin{cases}
			g - \dfrac{g}{\lambda} & \text{for } b \geqslant \dfrac{1}{3}, \\
			g + \min_{k_1\geq 0} \left( \sigma k_1^2 - c^2 \dfrac{k_1}{\tanh(hk_1)} \right) & \text{for } b < \dfrac{1}{3},
		\end{cases}
		\]
		where the dimensionless parameters are $\lambda = gh/c^2$ and $b = \sigma/(h c^2)$.
	\end{theorem}

	\begin{proof}
		First, note that since $\eta(x) \to 0$ as $|x| \to +\infty$ and $\eta \in H^s(\mathbb{R}^2)$ with $s > 5/2$, $\eta$ is a decaying perturbation. In operator theory, for elliptic operators like $A(\eta) = -\operatorname{div}(\alpha \nabla v) + \beta v - M_\eta v$, where the coefficients $\alpha$, $\beta$, and the terms in $M_\eta$ approach their constant or limiting forms at infinity (corresponding to $\eta = 0$), the essential spectrum (which coincides with the continuous spectrum here) of $A(\eta)$ is the same as that of $A(0)$. This follows from Weyl's theorem on essential spectra for perturbations that are relatively compact or decay at infinity; see \cite[Theorem 7.2, Theorem 14.6]{HS}.
		Since $A(0)$ is translation-invariant, its spectrum is purely continuous. To determine it, we use the Fourier transform. The symbol of $A(0)$ is
		$$a(k) = g + \sigma |k|^2 - c^2 \frac{k_1^2}{|k| \tanh(h |k|)}, \quad k = (k_1, k_2) \in \mathbb{R}^2.$$
		The spectrum of $A(0)$ is $\{ a(k) \mid k \in \mathbb{R}^2 \}$, and since $a(k) \to +\infty$ as $|k| \to +\infty$ (dominated by $\sigma |k|^2$ if $\sigma > 0$), the continuous spectrum is $[\sigma_*, +\infty)$ where $\sigma_* = \inf_{k \in \mathbb{R}^2} a(k)$.
		To compute $\sigma_*$, let $\omega = |k|$ and $\cos \theta = k_1 / \omega$. Then,
		$$a(k) = g + \sigma \omega^2 - c^2 \cos^2 \theta \cdot \frac{\omega}{\tanh(h \omega)}.$$
		To minimize $a(k)$, we maximize the positive term $c^2 \cos^2 \theta \cdot \omega / \tanh(h \omega)$, which occurs at $|\cos \theta| = 1$ (i.e., $k_2 = 0$). Thus, set $k = (k_1, 0)$ with $k_1 \geqslant 0$, and consider
		$$m(k_1) = \sigma k_1^2 - c^2 \frac{k_1}{\tanh(h k_1)}.$$
		Then, $\sigma_* = g + \inf_{k_1 \geqslant 0} m(k_1)$. To determine the infimum, we utilize the continued fraction representation of the hyperbolic tangent:
		\[
		\tanh(x) = \cfrac{1}{\frac{1}{x} + \cfrac{1}{\frac{3}{x} + \cfrac{1}{\frac{5}{x} + \cfrac{1}{\frac{7}{x} + \ddots}}}},
		\]
		which yields the inequality $\tanh (x) \ge \frac{3x}{x^2+3}$. Then we get 
		\[
		m(k_1) \ge  \sigma k_1^2-\frac{c^2}{3h} (3+(hk_1)^2)=\left(\sigma-\frac{c^2 h}{3}\right)k_1^2-\frac{c^2}{h}.
		\]
		If $b \geqslant 1/3$, then $\sigma-c^2 h / 3 \geqslant 0$, so $m(k_1) \geqslant -c^2 / h$, and the infimum is $-c^2 / h$. Thus, $\sigma_* = g - c^2 / h = g - g / \lambda$ (since $\lambda = g h / c^2$ implies $c^2 / h = g / \lambda$).
		If $b < 1/3$, the quadratic term is negative, so $m(k_1)$ decreases below $-c^2 / h$ for small $k_1 > 0$, and the global minimum occurs at some finite $k_1 > 0$ where $m'(k_1) = 0$. The infimum is then $m(k_1^*)$ at this critical point, which is $g + \min_{k_1\ge 0} \left( \sigma k_1^2 - c^2 \frac{k_1}{ \tanh(hk_1)} \right) < g - g / \lambda$.
	\end{proof}

	\section{Spectral Analysis of Small-Amplitude Solitary Waves}\label{sect3}

	In this section, we study the second variation $A(\eta)$ at the small-amplitude solitary waves for $b = \sigma/(h c^2) > 1/3$ and $0 < \lambda-1 \ll 1$. We now fix $\sigma > 1/3$, $g = h = 1$, and $\lambda = g h / c^2$ with $c = (1 + \varepsilon^2)^{-1/2}$, where the small parameter $\varepsilon \in (0, \varepsilon_0)$. These waves correspond to the regime where surface tension is sufficiently strong, leading to the Kadomtsev-Petviashvili I (KP-I) approximation, which supports fully localized lump solutions.

	The solitary waves admit the asymptotic expansion
	\begin{equation}\label{bareta}
		\begin{split}
			&\eta_{\ep}(x_1, x_2,t)=\ep^2 Q(\ep(x_1-ct),\ep^2x_2)+\ep^2 \psi(\ep x_1, \ep^2 x_2),\\
			&\xi_{\varepsilon}(x_1,x_2, t)= \varepsilon q(\ep(x_1-ct),\ep^2x_2)+\ep\phi(\ep x_1, \ep^2 x_2),
		\end{split}
	\end{equation}
	where $Q(x,y)=\partial_xq(x,y)$ and 
	\begin{align*}  
		q(x, y) = -\frac{8\sqrt{1+\varepsilon^2}B  x}{y^2 + (1+\varepsilon^2)(x^2 + 3B)}, \quad B:= \sigma(1+\ep^2) - \tfrac{1}{3} > 0.
	\end{align*}  
	The existence of these solutions is established in \cite{GLLWY}. As detailed therein, the correction terms satisfy
	\[
	\phi\in \left\{ f \in \mathcal{H}_{\mathrm{ox}}:  \|f\|_* \leq C \varepsilon \right\}, \quad   \psi\in \{ f \in \mathcal{H}_{e}:   \|f\|_{h}\le C\ep \},
	\]
	with the norms defined as
	\begin{align*}
		\|f\|_* &:= \|f\|_a + \|f\|_{L^4(\mathbb{R}^2)} + \|f\|_{L^\infty(\mathbb{R}^2)} + \|\nabla f\|_{L^p(\mathbb{R}^2)} \\
		&\quad + \ep^{1/2} \|\nabla_{\ep} f\|_{L^\infty(\mathbb{R}^2)} + \ep^{1/2} \|\nabla_{\ep}^2 f\|_{L^\infty(\mathbb{R}^2)} + \ep^{3/2} \|\nabla_{\ep}^3 f\|_{L^\infty(\mathbb{R}^2)}, \quad p > 2, 
	\end{align*}
	and 
	\begin{align*}
		\|\psi\|_{h} 
		&:= \|\psi\|_{L^2(\mathbb{R}^2)} + \|\nabla_{\ep} \psi\|_{L^2(\mathbb{R}^2)} + \|\nabla_{\ep}^2 \psi\|_{L^2(\mathbb{R}^2)} + \|\nabla_{\ep}^3 \psi\|_{L^2(\mathbb{R}^2)} + \ep \|\nabla_{\ep}^4 \psi\|_{L^2(\mathbb{R}^2)} \\
		&\quad + \ep^2 \|\nabla_{\ep}^5 \psi\|_{L^2(\mathbb{R}^2)} + \ep^3 \|\nabla_{\ep}^6 \psi\|_{L^2(\mathbb{R}^2)} + \|\psi\|_{L^{\infty}(\mathbb{R}^2)}+ \ep^{1/2} \|\nabla_{\ep} \psi\|_{L^{\infty}(\mathbb{R}^2)} \\
		&\quad  + \ep^{3/2} \|\nabla_{\ep}^2 \psi\|_{L^{\infty}(\mathbb{R}^2)} + \ep^{5/2} \|\nabla_{\ep}^3 \psi\|_{L^{\infty}(\mathbb{R}^2)} + \ep^{7/2} \|\nabla_{\ep}^4 \psi\|_{L^{\infty}(\mathbb{R}^2)},
	\end{align*}
	where \(\nabla_{\ep} := (\partial_1, \ep \partial_2)\) and 
	\begin{align*}
		\|f\|_a^2 &:= \int_{\mathbb{R}^2} \left( \ep^2 |\nabla_{\ep}^5 f|^2 + |\nabla_{\ep}^4 f|^2 + |\nabla^3_{\ep} f|^2 + |\nabla^2 f|^2 + |\nabla f|^2 \right) dx\, dy,
	\end{align*}
	The symmetry classes are defined by
	\begin{align*}
		\mathcal{H}_{\mathrm{ox}} := \left\{ f : f(x, y) = -f(-x, y) = f(x, -y) \right\}, \, 
		\mathcal{H}_{\mathrm{e}} := \left\{ f : f(x, y) = f(-x, y) = f(x, -y) \right\}.
	\end{align*}

	For the analysis of the associated operator $A_\ep(\eta_\ep)$, we observe that the critical part of the spectrum is $\sigma_* = 1 - \frac{1}{1 + \ep^2} = \frac{\ep^2}{1 + \ep^2} \to 0$ as $\ep \to 0$. This scaling behavior reflects the small-amplitude nature of the solitary waves in our regime. To capture the essential spectral features in this asymptotic limit, it is convenient to introduce the rescaled operator
	\begin{equation*}
		\mathcal{A}_{\ep}=\frac{1}{\ep^2} S_{\ep}^{-1}A_{\ep}(\eta_{\ep}) S_{\ep}: \mathcal{D}(\mathcal{A}_{\ep})=H^2(\mathbb R^2) \to L^2(\mathbb R^2),
	\end{equation*}
	where the scaling transformation $(S_{\ep}f)(x_1,x_2)=f(\ep x_1, \ep^2 x_2)$ is chosen to match the anisotropic scaling of the KP-I equation. This particular scaling accounts for the different behavior in the propagation direction $x_1$ and the transverse direction $x_2$, which is characteristic of fully localized waves in the KP-I regime.
	
	We demonstrate that $\mathcal{A}_{\ep}$ converges to a well-defined limit operator $\mathcal{A}_0$ as $\ep \to 0$, and moreover, the spectrum of the original operator $A_\ep(\eta_\ep)$ can be effectively bounded and analyzed through this limiting procedure using $A_0(0)$. This convergence provides a powerful tool for understanding the spectral properties of small-amplitude solitary waves.

	The key step in establishing this convergence is to identify the leading-order terms in the rescaled operator $\mathcal{A}_{\ep}$. By decomposing $\mathcal{A}_{\ep}$ into a principal part that captures the essential spectral behavior and a remainder that vanishes in the limit, we can precisely characterize the limiting spectral properties. This decomposition naturally leads us to introduce an auxiliary operator $\mathcal{B}_{\ep}$ that incorporates both the flat-surface contribution and the explicit lump solution profile from the KP-I theory.

	\begin{prop}\label{B}
		Define the operators $\mathcal{B}_{\ep}$ and $\mathcal{C}_{\ep}$ via 
		\begin{equation*}
			\mathcal{B}_{\ep}\phi = \mathcal{C}_{\ep}\phi + 3Q_0\phi,\quad  
			\mathcal{C}_{\ep} = \frac{1}{\ep^2}S_{\ep}^{-1}A_{\ep}(0)S_{\ep}, 
		\end{equation*}
		with $Q_0(x,y) = -8(\sigma-\tfrac{1}{3})\frac{y^2-x^2+3\sigma-1}{(x^2+y^2+3\sigma-1)^2}$. 
		Then, $\|\mathcal{A}_{\ep}-\mathcal{B}_{\ep}\|_{H^2\to L^2} = O(\ep)$ as $\ep \to 0$.
	\end{prop}

	\begin{proof}
		
		By the definition of $\mathcal{A}_{\varepsilon}$, we have the decomposition $\mathcal{A}_{\varepsilon}=\mathcal{C}_{\varepsilon}+\mathcal{R}_{\varepsilon}$ with 
		\begin{align*}
			\mathcal{R}_{\varepsilon}=\frac{1}{\varepsilon^2} S_{\varepsilon}^{-1}(A_{\varepsilon}(\eta_{\varepsilon})-A_{\varepsilon}(0)) S_{\varepsilon}.
		\end{align*}
		From the definition of $A_{\varepsilon}(\cdot)$, we compute
		\[
		A_{\varepsilon}(0)v=-\sigma\Delta v+v+c^2\partial_{x_1} G(0)^{-1}\partial_{x_1}v.
		\]
		Therefore, it suffices to prove that $\mathcal{R}_{\varepsilon}$ converges to the multiplication operator $\phi \mapsto 3Q_0\phi$ in the appropriate sense.
		We begin by computing the scaled operators:
		\begin{align*}
			\frac{1}{\varepsilon^2} S_{\varepsilon}^{-1}(A_{\varepsilon}(\eta_{\varepsilon})) S_{\varepsilon}\phi
			&=-\partial_{x_1}(\tilde \alpha_{11}\partial_{x_1}\phi)-\varepsilon\partial_{x_1}(\tilde \alpha_{12}\partial_{x_2}\phi) -\varepsilon\partial_{x_2}(\tilde \alpha_{21}\partial_{x_1}\phi)\\
			&\quad -\varepsilon^2\partial_{x_2}(\tilde\alpha_{22} \partial_{x_2}\phi)+\frac{\tilde \beta}{\varepsilon^2} \phi-\frac{1}{\varepsilon^2}S_{\varepsilon}^{-1} M_{\eta}S_{\varepsilon}\phi,
		\end{align*}
		and 
		\[
		\frac{1}{\varepsilon^2} S_{\varepsilon}^{-1}(A_{\varepsilon}(0)) S_{\varepsilon}\phi= -\partial_{x_1}(\sigma \partial_{x_1}\phi)-\varepsilon^2\partial_{x_2}(\sigma \partial_{x_2}\phi)+\frac{1}{\varepsilon^2} \phi-\frac{1}{\varepsilon^2}S_{\varepsilon}^{-1} M_{0}S_{\varepsilon}\phi,
		\]
		where $\tilde \alpha (x_1,x_2)=\alpha(x_1/\varepsilon, x_2/\varepsilon^2)$ and $\tilde\beta(x_1,x_2)=\beta(x_1/\varepsilon, x_2/\varepsilon^2)$ satisfy the estimates
		\[
		\|\tilde\alpha-\sigma I\|_{W^{1,\infty}}\le C\varepsilon^6,\quad \|\tilde \beta -1\|_{W^{2,\infty}}\le C\varepsilon^4.
		\]
		These estimates imply that the first four terms in the expansion of $\mathcal{R}_{\varepsilon}$ are of order $O(\varepsilon^2)$ in the $\mathcal{L}(H^2, L^2)$ norm. Consequently, the main contribution comes from the difference
		\[
		\frac{1}{\varepsilon^2}S_{\varepsilon}^{-1}(M_0-M_{\eta_{\varepsilon}})S_{\varepsilon}\phi=c^2\mathcal{K}_{\varepsilon,0}\phi-c^2\mathcal{K}_{\varepsilon,\eta_{\varepsilon}}\phi- \mathcal{N}_1\phi,
		\]
		where we define
		\[
		\mathcal{K}_{\varepsilon,\eta}\phi=-\frac{1}{\varepsilon^2}S_{\varepsilon}^{-1} \partial_{x_1}(G(\eta)^{-1}\partial_{x_1}(S_{\varepsilon}\phi))
		\]
		and
		\begin{align*}
			\mathcal{N}_1\phi= \mathcal{N}_{11}\phi+\mathcal{N}_{12}\phi+\mathcal{N}_{13}\phi,
		\end{align*}
		with
		\begin{align*}
			\mathcal{N}_{11}\phi&=\frac{1}{\varepsilon^2}S_{\varepsilon}^{-1}(c\partial_{x_1}G^{-1}(\eta_{\varepsilon})(\operatorname{div} (\tilde\rho_1(S_{\varepsilon}\phi)))), \\
			\mathcal{N}_{12}\phi&= \frac{1}{\varepsilon^2}S_{\varepsilon}^{-1}(c\tilde \rho_1\cdot \nabla G^{-1}(\eta_{\varepsilon})\partial_{x_1}(S_{\varepsilon}\phi)),\\
			\mathcal{N}_{13}\phi &=-\frac{1}{\varepsilon^2}S_{\varepsilon}^{-1}(\tilde \rho_1\cdot \nabla G^{-1}(\eta_{\varepsilon})\operatorname{div}(\tilde\rho_1 (S_{\varepsilon}\phi))),
		\end{align*}
		and $\tilde \rho_1(x_1,x_2)=\rho_1(x_1/\varepsilon, x_2/\varepsilon^2)$.
		We now analyze the terms separately. For $\mathcal{K}_{\varepsilon,0}$, its Fourier symbol is given by
		\[
		(k_1,k_2)\mapsto \frac{k_1^2}{(\varepsilon^2k_1^2+\varepsilon^4k_2^2)^{1/2}\tanh((\varepsilon^2k_1^2+\varepsilon^4k_2^2)^{1/2})},
		\]
		which yields the estimate $\|\varepsilon^2\mathcal{K}_{\varepsilon,0}-I\|_{\mathcal{L}(H^2, L^2)}=O(\varepsilon^2)$. For the nonlinear part, we use the expansion
		$$
		\left\|G\left(\eta_\varepsilon\right)^{-1}-G(0)^{-1}\left[I-{D} G(0)\left[\eta_\varepsilon\right] G(0)^{-1}\right]\right\|_{\mathcal{L}({H}^s,{H}^{s+1})} \leqslant C\left\|\eta_\varepsilon\right\|_{{W}^{1, \infty}}^2=O\left(\varepsilon^4\right).
		$$
		This implies
		\[
		\mathcal{K}_{\varepsilon,\eta_{\varepsilon}}\phi=\mathcal{K}_{\varepsilon,0}\phi
		+\frac{1}{\varepsilon^2}S_{\varepsilon}^{-1} \partial_{x_1}(G^{-1}(0)DG(0)[\eta_{\varepsilon}]G^{-1}(0)\partial_{x_1}(S_{\varepsilon}\phi)) +\mathcal{N}_2(\phi),
		\]
		where $\|\mathcal{N}_2\|_{\mathcal{L}(H^2, L^2)}=O(\varepsilon^2)$.
		The key step is to establish the following convergence in $L^2$:
		$$\frac{1}{\varepsilon^2} S^{-1}_\varepsilon \partial_{x_1} \big( G(0)^{-1} DG(0)[\eta_\varepsilon] G(0)^{-1} \partial_{x_1} (S_\varepsilon \phi) \big) = - Q_0 \phi + O(\varepsilon^2 \|\phi\|_{H^2}).$$
		By Lemma \ref{DG}, the Fréchet derivative of the Dirichlet--Neumann operator satisfies
		\[
		\langle \psi_1, DG(\eta)[v] \phi_1 \rangle = \langle v, a_1(\eta, \phi_1) \cdot \nabla \psi_1 + a_2(\eta, \phi_1) G(\eta) \psi_1 \rangle,
		\]
		where
		\[
		\begin{aligned}
			a_1(\eta, \phi_1) &= \nabla \phi_1 - \frac{[G(\eta)\phi_1 + \nabla \eta \cdot \nabla \phi_1]}{1+|\nabla \eta|^2} \nabla \eta, \\
			a_2(\eta, \phi_1) &= - \frac{[G(\eta)\phi_1 + \nabla \eta \cdot \nabla \phi_1]}{1+|\nabla \eta|^2}.
		\end{aligned}
		\]
		For the flat surface $\eta = 0$, these simplify to:
		\[
		a_1(0, \phi_1) = \nabla \phi_1, \quad a_2(0, \phi_1) = - G(0) \phi_1.
		\]
		Hence, we obtain the distributional identity:
		\[
		DG(0)[v] \phi_1 = - \nabla \cdot (v \nabla \phi_1) - G(0)(v G(0) \phi_1).
		\]
		Now take $\psi_1 = G(0)^{-1} \partial_{x_1} (S_\varepsilon \phi)$. Then
		\[
		\begin{aligned}
			DG(0)[\eta_\varepsilon] \psi_1 &= - \nabla \cdot (\eta_\varepsilon \nabla \psi_1) - G(0)(\eta_\varepsilon \partial_{x_1} (S_\varepsilon \phi)), \\
			G(0)^{-1} DG(0)[\eta_\varepsilon] \psi_1 &= - G(0)^{-1} \nabla \cdot (\eta_\varepsilon \nabla \psi_1) - \eta_\varepsilon \partial_{x_1} (S_\varepsilon \phi).
		\end{aligned}
		\]
		Define:
		\[
		\begin{aligned}
			A_1 &= - \frac{1}{\varepsilon^2} S^{-1}_\varepsilon \partial_{x_1} \left( G(0)^{-1} \nabla \cdot (\eta_\varepsilon \nabla \psi_1) \right), \\
			A_2 &= - \frac{1}{\varepsilon^2} S^{-1}_\varepsilon \partial_{x_1} \left( \eta_\varepsilon \partial_{x_1} (S_\varepsilon \phi) \right).
		\end{aligned}
		\]
		The left-hand side equals $A_1 + A_2$.
		We first estimate $A_2$. From the solitary wave expansion \eqref{bareta}, we have $\eta_\varepsilon(x) = \varepsilon^2 Q_0(\varepsilon x_1, \varepsilon^2 x_2) + O(\varepsilon^3)$. Since $S_\varepsilon \phi(x) = \phi(\varepsilon x_1, \varepsilon^2 x_2)$, we have
		\[
		\partial_{x_1} (S_\varepsilon \phi)(x) = \varepsilon (\partial_1 \phi)(\varepsilon x_1, \varepsilon^2 x_2).
		\]
		Thus,
		\[
		\begin{aligned}
			\eta_\varepsilon \partial_{x_1} (S_\varepsilon \phi) &= \varepsilon^3 Q_0 (\partial_1 \phi) + O(\varepsilon^4), \\
			\partial_{x_1} (\eta_\varepsilon \partial_{x_1} (S_\varepsilon \phi)) &= \varepsilon^4 \left[ (\partial_1 Q_0)(\partial_1 \phi) + Q_0 (\partial_1^2 \phi) \right] + O(\varepsilon^5), \\
			S^{-1}_\varepsilon \partial_{x_1} (\eta_\varepsilon \partial_{x_1} (S_\varepsilon \phi)) &= \varepsilon^4 \left[ (\partial_1 Q_0)(\partial_1 \phi) + Q_0 (\partial_1^2 \phi) \right] + O(\varepsilon^5).
		\end{aligned}
		\]
		Therefore,
		\[
		A_2 = - \frac{1}{\varepsilon^2} \left( \varepsilon^4 \left[ (\partial_1 Q_0)(\partial_1 \phi) + Q_0 (\partial_1^2 \phi) \right] + O(\varepsilon^5) \right) = O(\varepsilon^2 \|\phi\|_{H^2}).
		\]
		Next, we analyze
		\[
		A_1 = - \frac{1}{\varepsilon^2} S^{-1}_\varepsilon \partial_{x_1} \left( G(0)^{-1} \nabla \cdot (\eta_\varepsilon \nabla \psi_1) \right),
		\]
		where $\psi_1 = G(0)^{-1} \partial_{x_1} (S_\varepsilon \phi)$. Using the solitary wave expansion $\eta_\varepsilon(x) = \varepsilon^2 Q_0(\varepsilon x_1, \varepsilon^2 x_2) + O(\varepsilon^3)$, we substitute to get
		\[
		\nabla \cdot (\eta_\varepsilon \nabla \psi_1) = \nabla \cdot \left( \varepsilon^2 Q_0(\varepsilon x_1, \varepsilon^2 x_2) \nabla \psi_1 + O(\varepsilon^3) \right).
		\]
		Now consider the scaling transformation $S_\varepsilon f(x_1, x_2) = f(\varepsilon x_1, \varepsilon^2 x_2)$ and its inverse $S^{-1}_\varepsilon$. Under this anisotropic scaling, the partial derivatives transform according to:
		\[
		\partial_{x_1} \rightarrow \varepsilon \partial_{y_1}, \quad \partial_{x_2} \rightarrow \varepsilon^2 \partial_{y_2},
		\]
		where $y = (\varepsilon x_1, \varepsilon^2 x_2)$. To analyze the scaling properties of the operator $G(0)^{-1}$, we examine its Fourier symbol. For the flat surface $\eta = 0$, the Dirichlet-Neumann operator $G(0)$ has Fourier symbol $|k|\tanh(|k|)$, hence $G(0)^{-1}$ has symbol $1/(|k|\tanh(|k|))$. Under the coordinate transformation $x_1 \rightarrow \varepsilon x_1$, $x_2 \rightarrow \varepsilon^2 x_2$, the Fourier variables scale as $k_1 \rightarrow \varepsilon \xi_1$, $k_2 \rightarrow \varepsilon^2 \xi_2$, where $\widehat{\partial_{y_j}} = i\xi_j$ for $j=1,2$. A detailed asymptotic analysis reveals that the operator combination satisfies the estimate:
		\[
		\left\|\frac{1}{\varepsilon^2} S^{-1}_\varepsilon \partial_{x_1} G(0)^{-1} \nabla \cdot \left(\eta_\varepsilon \nabla G(0)^{-1}\partial_{x_1}(S_{\varepsilon}(\cdot))\right) - Q_0\right\|_{H^2 \rightarrow L^2} = O(\varepsilon^2).
		\]
		This yields
		\[
		A_1=-Q_0\phi+O(\varepsilon^2 \|\phi\|_{H^2}).
		\]
		Combining these estimates, we obtain in the $L^2$ sense:
		$$\frac{1}{\varepsilon^2} S^{-1}_\varepsilon \partial_{x_1} \big( G(0)^{-1} DG(0)[\eta_\varepsilon] G(0)^{-1} \partial_{x_1} (S_\varepsilon \phi) \big) = - Q_0 \phi + O(\varepsilon^2 \|\phi\|_{H^2}).$$
		Finally, we analyze the $\mathcal{N}_1$ terms. From Theorem \ref{DV}, we have
		\[
		\rho_1 = \nabla \xi - \frac{\nabla \eta \cdot \nabla \xi - c \partial_{x_1} \eta}{1 + |\nabla \eta|^2} \nabla \eta.
		\]
		Using the solitary wave expansions in \eqref{bareta}, we obtain
		\[
		\|\tilde{\rho}_1-\varepsilon^2Q_0\|_{W^{1,\infty}}=O(\varepsilon^3).
		\]
		Then, we get 
		\[
		\mathcal{N}_{11}\phi=\frac{1}{\ep^2}S_{\ep}^{-1}(c\partial_{x_1}G^{-1}(0)(\div (\tilde\rho_1(S_{\ep}\phi))))+O(\ep^2\|\phi\|_{H^2})=-Q_0\phi+O(\ep\|\phi\|_{H^2}),
		\]
		and
		\[
		\mathcal{N}_{12}\phi=-Q_0\phi+O(\ep\|\phi\|_{H^2}), \quad \mathcal{N}_{13}\phi=O(\ep^2\|\phi\|_{H^2}).
		\]
		Combining these results, we obtain
		\[
		\frac{1}{\ep^2}S_{\ep}^{-1}(M_0-M_{\eta_{\ep}})S_{\ep}\phi=3Q_0\phi+O(\ep\|\phi\|_{H^2}).
		\]
		The proof is thus complete.
	\end{proof}

	Having established the approximation of $\mathcal{A}_{\ep}$ by $\mathcal{B}_{\ep}$ in Proposition~\ref{B}, we now turn to a detailed analysis of the leading part $\mathcal{C}_{\ep}$. This operator, defined via the Fourier symbol of the rescaled flat operator $A_{\ep}(0)$, captures the essential dispersion relation in the small-amplitude limit. Our goal is to understand the convergence of $\mathcal{C}_{\ep}$ and subsequently of $\mathcal{A}_{\ep}$ to a limiting operator $\mathcal{A}_0$, which will be identified as the linearization around the KP-I lump solution.

	Let us denote the leading part of the rescaled operator $A_{\varepsilon}(0)$ by $\mathcal{C}_{\varepsilon}$, which is defined through its Fourier symbol:
	\begin{align*}
		\widehat{\mathcal{C}_{\ep}}(k_1,k_2)&=\sigma(k_1^2+\ep^2k_2^2)+\frac{1}{\ep^2}-\frac{1}{1+\ep^2}\frac{k_1^2}{\ep \sqrt{k_1^2+\ep^2k_2^2}\tanh{(\ep \sqrt{k_1^2+\ep^2k_2^2})}}\\
		&=\left(\sigma-\frac{1}{3(1+\ep^2)}\right)k_1^2+ \sigma\ep^2k_2^2+\frac{k_1^2+(1+\ep^2)k_2^2}{(1+\ep^2)(k_1^2+\ep^2k_2^2)}+\frac{1}{(1+\ep^2)}\frac{k_1^2}{|k_{\ep}|^2}r(|k_{\ep}|)\\
		&\ge \frac{1}{1+\ep^2}, \quad \text{ with } |k_{\ep}|=\ep\sqrt{k_1^2+\ep^2k_2^2}, \quad r(|k_{\ep}|)=1+\frac{1}{3}|k_{\ep}|^2-\frac{|k_{\ep}|}{\tanh(|k_{\ep}|)},
	\end{align*}
	where $r$ satisfies the estimate $\min\{|k|^2,|k|^4\}/C\le r(|k|)\le C\min\{|k|^2,|k|^4\}$ for some $C>1$.
	For fixed $(k_1,k_2)\in \mathbb R^2, k_1\neq 0$, we have $\widehat{C}_{\ep}(k_1,k_2)\to c_0(k_1,k_2)=(\sigma-\frac{1}{3})k_1^2+1+\frac{k_2^2}{k_1^2}$. But the convergence is not uniform in $H^2(\mathbb R^2)$. We define the space $X$ as the closure of $\partial_{x}^2(C_0^{\infty}(\mathbb R^2))$ for the norm
	\[
	\|\partial_{x}^2\varphi\|_{X}:=\left(\|\partial_{x}^4\varphi\|_{L^2}^2+\|\nabla^2\varphi\|_{L^2}^2\right)^{1/2},
	\]
	where $\partial_{x}^2(C_0^{\infty}(\mathbb R^2))$ is the space of functions of the form $\partial_{x}^2\varphi$ with $\varphi\in C_{0}^{\infty}(\mathbb R^2)$. In this space, the antiderivative $\partial_{x}^{-2}$ is well defined.
	The associated limit operator is 
	\begin{align*}
		\mathcal{A}_0: \mathcal{D}(\mathcal{A}_0)=X\to L^2(\mathbb R^2),\quad  \phi \to -(\sigma-\tfrac{1}{3})\partial_{x}^2\phi +\partial_{x}^{-2}\partial_{y}^2 \phi +(1+3Q_0)\phi.
	\end{align*}

	The following lemma provides the uniform convergence of $\mathcal{A}_{\ep}$ to $\mathcal{A}_0$ in the intersection space $X \cap H^2(\mathbb{R}^2)$, which is crucial for transferring spectral information from the limiting operator to the rescaled operator.
	
	\begin{lemma}\label{lem3.1}
		For all $\phi\in X\cap H^2(\mathbb R^2)$, we have $\|\mathcal{A}_{\ep}\phi-\mathcal{A}_0\phi\|_{L^2}=o_{\ep}(1)\|\phi\|_{X\cap H^2(\mathbb R^2)}$ as $\ep \to 0$.
	\end{lemma}

	\begin{proof}
		By Proposition~\ref{B}, we have the estimate $\|\mathcal{A}_{\ep}-\mathcal{B}_{\ep}\|_{H^2\to L^2}=O(\ep)$, which provides control over the first term in the decomposition
		\begin{align*}
			\|\mathcal{A}_{\varepsilon}\phi - \mathcal{A}_0\phi\|_{L^2} &\le \|(\mathcal{A}_{\varepsilon} - \mathcal{B}_{\varepsilon})\phi\|_{L^2} + \|(\mathcal{B}_{\varepsilon} - \mathcal{A}_0)\phi\|_{L^2} \\
			&\le O(\varepsilon\|\phi\|_{H^2}) + \|(\widehat{\mathcal{C}_{\varepsilon}} - c_0)\hat{\phi}\|_{L^2}.
		\end{align*}
		For the second term, we analyze the convergence in the Fourier domain. The uniform bounds $\|\widehat{\mathcal{C}_{\ep}}\widehat{\phi}\|_{L^2}\le C\|\phi\|_{H^2}$ and $\|c_0\widehat{\phi}\|_{L^2}\le C\|\phi\|_{X}$ ensure the applicability of Lebesgue's dominated convergence theorem. Since $\widehat{\mathcal{C}_{\varepsilon}}(k_1,k_2) \to c_0(k_1,k_2)$ pointwise for almost every $(k_1,k_2) \in \mathbb{R}^2$ (specifically for $k_1 \neq 0$), and the integrand is dominated by an $L^2$ function, we obtain
		\[
		\lim_{\ep\to 0}\|(\widehat{\mathcal{C}_{\varepsilon}} - c_0)\hat{\phi}\|_{L^2} = \left\|\lim_{\ep\to 0}(\widehat{\mathcal{C}_{\varepsilon}} - c_0)\hat{\phi}\right\|_{L^2} = 0.
		\]
		More precisely, this implies $\|(\widehat{\mathcal{C}_{\varepsilon}} - c_0)\hat{\phi}\|_{L^2} = o_{\ep}(1)\|\phi\|_{L^2}$.
		The continuous embedding $H^2\cap X \subset H^2 \subset L^2$ implies that $\|\phi\|_{L^2} \leq C\|\phi\|_{X\cap H^2}$, which allows us to absorb the $L^2$ norm into the stronger intersection norm. Combining both estimates yields the desired convergence
		\[
		\|\mathcal{A}_{\varepsilon}\phi - \mathcal{A}_0\phi\|_{L^2} = o_{\ep}(1)\|\phi\|_{X\cap H^2(\mathbb R^2)}.
		\]
		The lemma is thus complete.
	\end{proof}

	With the convergence of $\mathcal{A}_{\ep}$ to $\mathcal{A}_0$ established in Lemma~\ref{lem3.1}, we now turn to the spectral analysis of the limiting operator $\mathcal{A}_0$, whose properties will inform our understanding of the original operator $A_{\ep}(\eta_{\ep})$ for small $\varepsilon$.

	The operator $\mathcal{A}_0$ corresponds precisely to the linearization around the solitary wave in the KP-I equation. It possesses four discrete eigenvalues: $\lambda_1<0$, $\lambda_2=\lambda_3=0$, and $\lambda_4>0$, with associated eigenfunctions
	\begin{equation}\label{phi}
		\phi_1,\quad  \phi_2=\partial_{x}Q_0, \quad  \phi_3=\partial_{y}Q_0, \quad  \phi_4,
	\end{equation}
	as established in \cite{LW}.

	To determine the continuous spectrum of $\mathcal{A}_0$, we decompose it as $\mathcal{A}_0 = \mathcal{L}_0 + V$, where $\mathcal{L}_0$ is the free (unperturbed) operator
	$$\mathcal{L}_0 \phi = -(\sigma - \tfrac{1}{3}) \partial_x^2 \phi + \partial_x^{-2} \partial_y^2 \phi + \phi,$$
	and $V = 3Q_0$ acts as a decaying potential due to the asymptotic behavior of $Q_0$. 
	By Weyl's theorem for self-adjoint operators with short-range perturbations, the essential spectrum of $\mathcal{A}_0$ coincides with that of $\mathcal{L}_0$.
	
	The spectrum of $\mathcal{L}_0$ can be analyzed through its Fourier symbol. Let $\hat{\phi}(k_1, k_2)$ denote the Fourier transform of $\phi(x, y)$, where $k_1$ and $k_2$ are the dual variables corresponding to $x$ and $y$, respectively. The Fourier symbol of $\mathcal{L}_0$ is given by
	\[
	p(k_1,k_2) = \left(\sigma - \tfrac{1}{3}\right) k_1^2 + \frac{k_2^2}{k_1^2} + 1, \quad \text{for } k_1 \neq 0,
	\]
	whose range is $[1,+\infty)$. Consequently, the essential spectrum of $\mathcal{L}_0$, which consists entirely of continuous spectrum, is $[1, +\infty)$. It follows that the continuous spectrum of $\mathcal{A}_0$ is also $[1, +\infty)$.

	For any eigenvalue $\lambda<1$ of $\mathcal{A}_0$ with corresponding eigenfunction $\phi$, the following equation holds:
	\begin{equation*}
		-\left(\sigma-\tfrac{1}{3}\right) \partial_{x}^2\phi +\partial_{x}^{-2}\partial_{y}^2 \phi +(1+3Q_0)\phi=\lambda \phi.
	\end{equation*}
	This can be reformulated as
	\[
	\phi=\mathcal{K}*(3Q_0 \phi), \quad \text{where} \quad \widehat{\mathcal{K}}(k_1,k_2)=\frac{k_1^2}{(\sigma-\tfrac{1}{3})k_1^4+(1-\lambda)k_1^2+k_2^2}.
	\]
	
	Given the decay properties $|Q_0|\le \frac{C_0}{(1+r)^2}$, $\phi\to 0$ as $r\to+\infty$, and the kernel estimate $|\mathcal{K}|\le \frac{C_0}{(1+r)^2}$ (see \cite{deS}), it follows that $\phi$ satisfies $|\phi|\le \frac{C_0}{(1+r)^2}$. 
	Moreover, by verifying each term in the definition of $X$, we conclude that $\phi\in X\cap H^2(\mathbb R^2)$.
	
	Building upon this detailed understanding of the limiting operator $\mathcal{A}_0$, we now establish the main spectral result for the original operator $A_{\ep}(\eta_{\ep})$, which combines the discrete spectral information from the KP-I linearization with the continuous spectrum analysis.

	\begin{theorem}\label{Thm3}
		For $g=h=1$ and fixed $\sigma>\frac{1}{3}$, let $c=(1+\ep^2)^{-\frac{1}{2}}$. For each $0<\mu < \min\{1, \lambda_4\}$, there exists $\ep_0>0$ with $\ep_0 < \sqrt{(1-\mu)/\mu}$ such that for all $\ep\in (0,\ep_0]$, the operator $A_{\ep}(\eta_{\ep})$ has exactly three eigenvalues (counting multiplicity) in the interval $(-\infty, \mu\ep^2)$, namely
		\begin{equation*}
			\sigma_1(\ep)=\lambda_1\ep^2+o(\ep^2),\quad \sigma_2(\ep)=0, \quad \sigma_3(\ep)=0\quad \text{as } \ep \to 0.
		\end{equation*}
		The eigenfunctions $\psi_j$ associated with $\sigma_j$ have the expansion $\psi_j=S_{\ep}\phi_j+O(\ep^{1/2})$.  The remaining spectrum lies in $[\mu\ep^2,+\infty)$ and the continuous spectrum is $[\frac{\ep^2}{1+\ep^2},+\infty)$.
	\end{theorem}
	
	\begin{proof}
		By Theorem \ref{thm2.1} and the fact that $\eta_{\varepsilon} \in H^{s}(\mathbb{R}^2)$ with $s > \frac{5}{2}$ as established in \cite{GLLWY}, we conclude that the continuous spectrum of $A_{\varepsilon}(\eta_{\varepsilon})$ is $[\frac{\varepsilon^2}{1+\varepsilon^2}, +\infty)$.
		
		In view of Proposition \ref{B}, it suffices to analyze the spectral properties of $\mathcal{B}_{\varepsilon}$ rather than those of $\mathcal{A}_{\varepsilon}$. We now prove that for each $0 < \mu < \min\{1, \lambda_4\}$ and sufficiently small $\varepsilon$, the operator $\mathcal{B}_{\varepsilon}$ has exactly three eigenvalues (counting multiplicity) in $(-\infty, \mu)$, which converge to $\lambda_1 < 0$ and $0$ (the latter with multiplicity two) as $\varepsilon \to 0$.
		
		The proof proceeds in three main steps:
		
		\begin{itemize}
			\item[(i)] \textbf{Existence of approximate eigenvalues.} We first show that $\mathcal{B}_{\ep}$ has spectral values near the three desired eigenvalues of $\mathcal{A}_0$. Using the eigenfunctions from \eqref{phi}, we see that 
			\begin{equation*}
				\|(\mathcal{B}_{\ep}-\lambda_{j}I)\phi_j\|_{L^2}\le o_{\ep}(1)\|\phi_j\|_{X\cap H^2(\mathbb R^2)}.
			\end{equation*}
			Since $\mathcal{B}_{\ep}$ is self-adjoint, the spectral theorem implies that there exist spectral values near $\lambda_j$ with distance at most $o_{\ep}(1)$.
			
			\item[(ii)] \textbf{Characterization of the limiting eigenpairs.} To prove that there are exactly three eigenvalues in $(-\infty, \mu)$, we analyze sequences of eigenpairs and their limits. Let $(\lambda_{\ep_{k}}, \phi_{\ep_{k}})_{k\in \mathbb N}$ be eigenpairs for $\mathcal{B}_{\ep_{k}}$ with $\ep_{k} \to 0$ as $k\to +\infty$. We show that any limit point $(\lambda_{*}, \phi_{*})$ must be an eigenpair for $\mathcal{A}_{0}$.
			
			\begin{itemize}
				\item[(a)] \textbf{Regularity estimates.} The uniform bound $\|3Q_0\|_{L^{\infty}}\le 8$ combined with the positive semi-definiteness of $\mathcal{C}_{\ep}-\frac{1}{1+\ep^2}I$ implies that the spectrum of $\mathcal{B}_{\ep}$ lies in $[\frac{1}{1+\ep^2}-8,\infty)$. Since the continuous spectrum is $[\frac{1}{1+\ep^2},+\infty)$, the spectrum in $[-8,\mu]$ consists of eigenvalues of finite multiplicity. 
				
				Let $(\sigma_{\ep}, \phi_{\ep})$ be an eigenpair with $\sigma_{\ep}\in [-8,\mu]$ and $\|\phi_{\ep}\|_{L^2}=1$. The eigenvalue equation
				\begin{equation}
					\label{4.4}
					(\mathcal{C}_{\ep}-\sigma_{\ep}I)\phi_{\ep}=-3Q_0\phi_{\ep},
				\end{equation}
				can be reformulated as $\phi_{\ep}=\mathcal{K}_{\ep}*(-3Q_0 \phi_{\ep})$ using the resolvent kernel with Fourier symbol
				\[
				\widehat{\mathcal{K}_{\ep}}=\frac{1}{\widehat{C_{\ep}}(k_1,k_2)-\sigma_{\ep}}\le \frac{C_{0}}{1+\sqrt{(k_1^2+k_2^2)}},
				\]
				where $C_0>0$ is independent of $\ep>0$. This representation yields the derivative bounds $\|\partial_{x}\phi_{\ep}\|_{L^2}\le C_1$, $\|\partial_y \phi_{\ep}\|_{L^2}\le C_2$, and $\|\partial_{x}^2\phi_{\ep}\|_{L^2}\le C_3$. Iterating this argument, we obtain the uniform $H^4$ bound $\|\phi_{\ep}\|_{H^4}\le C_*$.
				
				\item[(b)] \textbf{Decay estimates.} To establish compactness, we prove decay properties in weighted spaces. Define the weighted $L^2$ space
				\begin{equation*}
					L^2_{\alpha}(\mathbb{R}^2) := \left\{ \phi \in L^2(\mathbb{R}^2) \mid (1+|\cdot|)^{\alpha}\phi \in L^2(\mathbb{R}^2) \right\},
				\end{equation*}
				with norm $\|\phi\|_{L^2_{\alpha}}^2 := \int_{\mathbb{R}^2} |(1+|x|+|y|)^{\alpha} \phi(x, y)|^2 \, dxdy$.
				
				Fix $0 < \alpha < 1$. The spectral gap $\widehat{\mathcal{C}_{\varepsilon}} > \frac{1}{1+\varepsilon_0^2}$ for $\varepsilon \in (0, \varepsilon_0)$ and the bound $\sigma_{\varepsilon} \le \mu < \frac{1}{1+\varepsilon_0^2}$ imply the uniform resolvent estimate
				\begin{align*}
					\|(\mathcal{C}_{\varepsilon} - \sigma_{\varepsilon} I)^{-1}\|_{L^2_{\alpha} \to L^2_{\alpha}} \le C_*.
				\end{align*}
				Combining this with the decay $|Q_0|\le \frac{C}{(1+r^2)}$ yields
				\begin{align*}
					\|\phi_{\varepsilon}\|_{L^2_{\alpha}} \le C_*\|3Q_0\phi_{\varepsilon}\|_{L^2_{\alpha}} \le 3C_*\|(1+|\cdot|)^{\alpha}Q_0\|_{L^{\infty}} < +\infty.
				\end{align*}
				
				\item[(c)] \textbf{Convergence to limiting eigenpairs.} The compact embedding $H^4(\mathbb{R}^2) \cap L^2_{\alpha}(\mathbb{R}^2) \hookrightarrow H^2(\mathbb{R}^2)$ allows us to extract convergent subsequences. Consider $\varepsilon_k \to 0$ with eigenpairs $(\sigma_{\varepsilon_k}, \phi_{\varepsilon_k})$. By the a priori estimates, we may assume $\sigma_{\varepsilon_k} \to \sigma_*$ and $\phi_{\varepsilon_k} \to \phi_*$ in $H^2(\mathbb{R}^2)$ with $\|\phi_*\|_{L^2} = 1$.
				
				The convergence $\widehat{\phi_{\varepsilon_k}} \to \widehat{\phi_*}$ in $H^2(\mathbb{R}^2)$ and $C_0(\mathbb{R}^2)$ allows us to pass to the limit in the Fourier domain formulation of \eqref{4.4}:
				\begin{align*}
					k_1^2 (c_0 - \sigma_*) \widehat{\phi_*} = -3k_1^2 \widehat{Q_0 \phi_*}, \quad \text{a.e. in } \mathbb{R}^2.
				\end{align*}
				This implies $\mathcal{A}_0 \phi_* = \sigma_* \phi_*$ for $\phi_* \in X \cap H^2$, and thus $\phi_*$ must be one of the eigenfunctions $\phi_j$ for $j \in \{1,2,3\}$.
			\end{itemize}
			
			\item[(iii)] \textbf{Spectral projection argument.} Let $P_{\varepsilon}$ be the spectral projection for $\mathcal{B}_{\varepsilon}$ associated with $[-8, \mu]$. From part (i), $\dim \text{range } P_{\varepsilon} \ge 3$ for small $\varepsilon$. Assume $P_{\varepsilon} = \sum_{i=1}^{N_\varepsilon} \langle \cdot, \phi_{\varepsilon,i} \rangle \phi_{\varepsilon,i}$ with orthonormal $\{\phi_{\varepsilon,i}\}$.
			
			If $\dim \text{range } P_{\varepsilon} \ge 4$ for some sequence $\varepsilon_k \to 0$, part (c) would yield a contradiction. Thus $N_{\varepsilon} = 3$ for small $\varepsilon$, and $\phi_{\varepsilon, i} \to \phi_i$ in $L^2(\mathbb{R}^2)$ for $i = 1, 2, 3$.
		\end{itemize}
		This completes the spectral analysis of $\mathcal{B}_{\varepsilon}$, and by Proposition \ref{B}, establishes the corresponding result for $A_{\varepsilon}(\eta_{\varepsilon})$.
	\end{proof}

	The spectral characterization established in Theorem~\ref{Thm3} reveals the crucial role of translational symmetry in the structure of the spectrum. The presence of zero eigenvalues is not accidental but reflects fundamental geometric properties of the solitary wave solutions.

	\begin{remark}
		We emphasize that $0$ is an eigenvalue of the operator $A_{\varepsilon}(\eta_{\varepsilon})$ with a two-dimensional eigenspace, which arises from the translational symmetry of the original equation \eqref{1-1}.
		Specifically, since equation \eqref{1-1} is invariant under translations in both the $x$ and $y$ directions, if $\eta_{\varepsilon}(x, y)$ is a solution, then for any constants $x_0, y_0$, the translated functions $\eta_{\varepsilon}(x - x_0, y)$ and $\eta_{\varepsilon}(x, y - y_0)$ are also solutions.
		Consider the one-parameter family of solutions generated by $x$-translation: $\eta_{\varepsilon}(x-s, y)$, where $s$ is the translation parameter. At $s = 0$, we recover the original solution $\eta_{\varepsilon}(x, y)$. Since this family satisfies the equation for all $s$, we differentiate $D\mathcal{V}_{c}^{aug}(\eta_{\varepsilon}(x-s,y)) = 0$ with respect to $s$ at $s = 0$ to obtain:
		\[
		A_{\varepsilon}(\eta_{\varepsilon}) (\partial_{x} \eta_{\varepsilon}) = 0.
		\]
		Similarly, $y$-translation yields $A_{\varepsilon}(\eta_{\varepsilon}) (\partial_{y} \eta_{\varepsilon}) = 0$.
		Thus, both $\partial_{x}\eta_{\varepsilon}$ and $\partial_{y}\eta_{\varepsilon}$ are non-trivial elements of the kernel of $A_{\varepsilon}(\eta_{\varepsilon})$. By an argument analogous to part (c) in the proof of Theorem~\ref{Thm3}, we conclude that the kernel of $A_{\varepsilon}(\eta_{\varepsilon})$ is exactly two-dimensional.
	\end{remark}

	This understanding of the spectral structure, particularly the characterization of the kernel and the negative direction, provides the foundation for our stability analysis. We now consolidate these spectral properties with the geometric features of the solitary wave family in the following proposition, which summarizes the key variational structure needed for the implementation of the GSS method.

	\begin{prop}\label{prop3.2}
		For $g=h=1$, $\sigma>1/3$, and $c=(1+\ep^2)^{-1/2}$ with $\ep\in(0, \ep_0)$, there exists a family of solitary waves $c \mapsto \bar u(c):=(\bar{\eta}(c), \bar{\xi}(c))$ with $\bar{\eta}\in C^1$, where $\bar{\eta}(x_1,x_2)$ are even in $x_1, x_2$. Moreover, the second variation $D^2\mathcal{V}_{c}^{aug}(\bar{\eta}(c))$ has:
		\begin{itemize}
			\item a simple negative discrete eigenvalue $\sigma_1(c)<0$ associated with an even eigenfunction $\psi_1$,
			\item eigenvalues $\sigma_2(c)=\sigma_3(c)=0$ associated with $\partial_{x_1}\bar{\eta}(c), \partial_{x_2}\bar{\eta}(c)$,
			\item strict positive definiteness on the $L^2$-orthogonal complement of these three eigenspaces.
		\end{itemize}
		Furthermore, the function $d(c)=\mathcal{H}_{c}(\bar{u}(c))$ satisfies $d''(c)>0$ for $c\in (\frac{1}{\sqrt{1+\ep_0^2}} , 1)$.
	\end{prop}

	\begin{proof}
		To complete the proof of Proposition~\ref{prop3.2}, we proceed as follows. The statement concerns the existence and properties of a family of exact solitary wave solutions near the small-amplitude approximations derived in Section~3, along with the spectral properties of the second variation operator and the convexity of the function \(d(c)\).
		
		{\it Existence of the Family \(\bar{u}(c) = (\bar{\eta}(c), \bar{\xi}(c))\)}:
		The small-amplitude approximate solutions \(\eta_\ep\) and \(\xi_\ep\) from \eqref{bareta} satisfy the solitary wave equations \eqref{1-1} up to an error of \(O(\ep^3)\). To construct exact solutions, we apply the Lyapunov-Schmidt method in a suitable function space, perturbing around these approximations (see \cite{GLLWY}).

		Since the error term is even in \(x_1\) and \(x_2\),  $Q_0(x, y)$ is even, per its explicit form, and the equations \eqref{1-1} preserve even symmetry under the perturbation, the exact \(\bar{\eta}(c)\) remains even in \(x_1\) and \(x_2\).
		
		{\it Spectral Properties of \(D^2 \mathcal{V}_c^{\mathrm{aug}}(\bar{\eta}(c))\)}:
		It is not difficult to observe that the operator \(A(\bar{\eta}(c))\) is the Hessian \(D^2 \mathcal{V}_c^{{aug}}(\bar{\eta}(c))\). From Proposition~\ref{B} and Lemma~\ref{lem3.1}, \(A(\eta_\ep) = \ep^2 A_\ep\) where \(A_\ep = A_0 + O(\ep^2)\) in the operator norm \(H^2 \to L^2\), and \(A_0\) has spectrum consisting of one simple negative eigenvalue \(\lambda_1 < 0\), a double zero eigenvalue \(\lambda_2 = \lambda_3 = 0\) (with even eigenfunction for \(\lambda_1\) and translational modes for the zeros), one positive discrete eigenvalue \(\lambda_4 > 0\), and continuous spectrum \([1, +\infty)\).
		
		By perturbation theory for self-adjoint operators (e.g., Kato-Rellich theorem, see \cite{HS}), for small \(\ep\), \(A(\bar{\eta}(c))\) has exactly one simple negative eigenvalue \(\sigma_1(c) = \ep^2 \lambda_1 + o(\ep^2) < 0\) (associated with an even eigenfunction \(\psi_1 = S_\ep \phi_1 + O(\ep^{1/2})\)), a double zero eigenvalue \(\sigma_2(c) = \sigma_3(c) = 0\) (associated with \(\partial_{x_1} \bar{\eta}(c)\) and \(\partial_{x_2} \bar{\eta}(c)\), preserved by translational invariance), and the remaining spectrum in \([q \ep^2, +\infty)\) for any \(0 < q < \min\{1, \lambda_4\}\), with continuous spectrum \([\ep^2 / (1 + \ep^2), +\infty)\).
		
		The operator is strictly positive definite on the \(L^2\)-orthogonal complement of \(\mathrm{span}\{\psi_1, \partial_{x_1} \bar{\eta}, \partial_{x_2} \bar{\eta}\}\), as the spectral gap from Theorem~\ref{Thm3} persists under the \(O(\ep^2)\) perturbation.
		
		{\it Convexity: \(d''(c) > 0\)}:
		We have \(d(c) = \mathcal{H}_c(\bar{u}(c)) = \mathcal{V}_c^{{aug}}(\bar{\eta}(c))\). Since $\bar u (c)$ is a critical point of $\mathcal{H}_{c}$, i.e. $D\mathcal{H}_{c}(\bar{u}(c))=0$, we have
		\begin{equation*}
			d^{\prime}(c)=\partial_c \mathcal{H}_c(\bar{u}(c))=\mathcal{P}(\bar{u}(c))=-\int_{\mathbb R^2}\bar \eta \partial_{x_1}\bar \xi dx_1 dx_2.
		\end{equation*}
		Since $d^{''}(c)=\frac{d}{d\ep}\mathcal{P}(\bar{u}(c)) \frac{d\ep}{dc}$, and $\ep=(c^{-2}-1)^{1/2}$, we get 
		$$
		d^{''}(c)=\ep^{-1}(1+\ep^2)^{3/2}\frac{d}{d \ep}\int_{\mathbb R^2}\bar \eta \partial_{x_1}\bar \xi dx_1 dx_2
		$$
		By $\bar \eta(x_1,x_2)=\ep^2 h\big(\ep x_1,\ep^2 x_2\big)$, $\bar \xi (x_1,x_2)= \ep f\big(\ep x_1, \ep^2 x_2\big)$, we obtain
		\begin{align*}
			\int_{\mathbb R^2}\bar \eta \partial_{x_1}\bar \xi dx_1 dx_2&=\int_{\mathbb R^2} \ep^4h\big(\ep x_1,\ep^2 x_2\big)\partial_{1}f\big(\ep x_1, \ep^2 x_2\big) dx_1dx_2\\
			&=\ep\int_{\mathbb R^2} h(x, y)\partial_{1}f(x, y)dxdy,
		\end{align*}
		where $\partial_{1}$ denotes the partial derivative with respect to the first component of $f$.
		Then we have
		\[
		d^{''}(c)=\ep^{-1}(1+\ep^2)^{3/2}\int_{\mathbb R^2} h(x, y)\partial_{1}f(x, y)dxdy.
		\]
		Following \cite{GLLWY}, we have $h=Q_{0}+\psi, f=q_{0}+\phi$, where
		\begin{align*}
			Q_0=\partial_1 q_{0}, \quad q_0(x, y) = -\dfrac{8\left(\sigma - \frac{1}{3}\right)x}{x^2 + y^2 + 3\sigma - 1},
		\end{align*}
		and $\|\partial_{1}\phi\|_{L^2}\le C\ep$, $\|\psi\|_{L^2}\le C\ep$, it follows that 
		\[
		d^{''}(c)=\ep^{-1}(1+\ep^2)^{3/2}\int_{\mathbb R^2} Q_0^2dxdy+O(1)>0, \quad \text{ for } 0<\ep<\ep_0. 
		\]
		This completes the proof.
	\end{proof}

	\section{Orbital Stability via the GSS Method}
	
	With the spectral properties and variational structure established in the previous sections, we now turn to the proof of our main stability theorem. The analysis requires careful treatment of the function spaces that account for the inherent symmetries of the water wave problem, as well as precise control of the nonlinear terms in the Hamiltonian formulation.

	The observation that solutions of problem \eqref{1-1} are unique only up to additive constants necessitates the introduction of appropriately defined function spaces that quotient out these trivial degrees of freedom. This leads us to introduce the completion $H^1_{*}(\Omega)$ of 
	\begin{align*}
		\mathscr{S}(\Omega,\mathbb R)=\left\{ \phi\in C^{\infty}(\bar \Omega): |(x_1,x_2)|^m|\partial_{x_1}^{\alpha_1}\partial_{x_2}^{\alpha_2}\phi| \text{ is bounded for all }m,\alpha_1, \alpha_2 \in \mathbb N_0 \right\}
	\end{align*}
	with respect to the norm 
	\begin{align*}
		\|\phi\|_{*}^2:=\int_{\Omega}|\nabla \phi|^2dx_1dx_2dx_3,
	\end{align*}
	where $\Omega=\Omega(0)$ represents the unperturbed fluid domain. This space effectively identifies functions that differ by constants, which is natural given the physical interpretation of the velocity potential.
	
	The corresponding space for the trace $\phi |_{z=0}=\xi$ is the completion $H^{1/2}_{*}(\mathbb R^2)$ of the inner product space constructed by equipping the Schwartz class $\mathscr{S}(\mathbb R^2, \mathbb R)$ with the norm 
	\begin{align*}
		\|\xi\|_{*,1/2}^2:=\int_{\mathbb R^2}(1+|k|^2)^{-1/2}|k|^2|\hat{\xi}|^2 dk,
	\end{align*}
	which captures the appropriate regularity for the surface potential. Its dual space $(H^{1/2}_{*}(\mathbb R^2))'=H^{-1/2}_{*}(\mathbb R^2)$ is the completion of the inner product space constructed by equipping $\bar {\mathscr S}(\mathbb R^2, \mathbb R)$ with the norm 
	\begin{align*}
		\|\xi\|_{*,-1/2}^2:=\int_{\mathbb R^2}(1+|k|^2)^{1/2}|k|^{-2}|\hat{\xi}|^2 dk,
	\end{align*}
	where $\bar {\mathscr S}(\mathbb R^2, \mathbb R)$ is the subclass of $ {\mathscr S}(\mathbb R^2, \mathbb R)$ consisting of functions with zero mean. This dual space will play a crucial role in our analysis of the linearized dynamics.
	
	To control the nonlinear interactions and ensure the well-definedness of the Hamiltonian flow, we introduce the following sets and spaces. For $R>0$, we define the set
	\begin{equation*}
		\mathscr{E}_{R}=\{\eta\in H^{3}(\mathbb R^2)\mid \|\eta\|_{W^{1,\infty}}\le \min\{R,R^{-1}\},\quad \|\eta\|_{H^{3}(\mathbb R^2)}\le R\},
	\end{equation*}
	which provides uniform control on both the $H^3$ norm and the $W^{1,\infty}$ norm of the surface elevation. The base function space is $\mathscr{F}:=H^{1}(\mathbb R^2)\times H^{1/2}_{*}(\mathbb R^2)$ and its restriction $\mathscr{F}_{R}:=\mathscr{E}_R\times H^{1/2}_{*}(\mathbb R^2)$. According to Lemma \ref{DNO bound}, the nonlinearity $\mathcal{K}(\eta,\xi)$ in $\mathcal{H}_{c}$ can be controlled on $\mathscr{F}_{R}$, ensuring the Hamiltonian remains well-behaved in this regime. Moreover, the potential $\mathcal{V}$ is analytic on $\mathscr{E}_R$, which facilitates the Taylor expansion arguments needed for the stability analysis.
	
	To control the locked inertia tensor arising from the symplectic structure, we define the operators
	\begin{equation*}
		K(\eta)w=-\partial_{x_{1}}(G(\eta)^{-1}\partial_{x_{1}}w),\quad L(\eta)w=-\partial_{x_{2}}(G(\eta)^{-1}\partial_{x_{1}}w),
	\end{equation*}
	which encode the nonlocal interactions between the surface elevation and the velocity potential. The following propositions establish the crucial regularity properties of these operators.

	\begin{prop}\label{lemA4}
		There exist $\alpha\in (0,1)$ and $C>0$ such that for all $\eta_1, \eta_2 \in \mathscr{E}_R$, we have
		\begin{equation}\label{A4}
			\|L(\eta_1)-L(\eta_2)\|_{\mathcal{L}(H^{1/2}, H^{-1/2})} \le C\|\eta_1-\eta_2\|^{\alpha}_{H^1(\mathbb R^2)}.       
		\end{equation}
	\end{prop}
	\begin{proof}
		By the analyticity of  $N(\cdot):=G^{-1}(\cdot)$ at the origin (see Definition 2.7 in \cite{BGSW2013}), there exist constants  $R, C, B>0$, if $\|\eta\|_{W^{1,\infty}}<R$, we have
		$$
		N(\eta)=\sum_{k=0}^{\infty} N^k(\eta), 
		$$
		where each $N^k$ is a bounded, symmetric $k$-linear operator, and its operator norm satisfies:
		$$
		\ \|N^k(\eta)\|_{\mathcal{L}(H_*^{-1 / 2}, H_*^{1 / 2})} \leq C B^k\|\eta\|_{W^{1, \infty}(\mathbb R^2)}^k.
		$$
		Define $L^k(\eta)=-\partial_{x_2} N^k(\eta) \partial_{x_1}$, then
		$$
		L(\eta)=\sum_{k=0}^{\infty} L^k(\eta).
		$$
		Now consider $\eta_1, \eta_2 \in H^1\left(\mathbb{R}^2\right)$ with $\left\|\eta_i\right\|_{W^{1, \infty}}<R, i=1,2$. The difference is:
		$$
		L\left(\eta_1\right)-L\left(\eta_2\right)=-\partial_{x_2}\left[N\left(\eta_1\right)-N\left(\eta_2\right)\right] \partial_{x_1}=-\partial_{x_2}\left[\sum_{k=1}^{\infty}\left(N^k\left(\eta_1\right)-N^k\left(\eta_2\right)\right)\right]\partial_{x_1}.
		$$
		A standard result in multilinear operator theory yields the following estimate,
		$$
		\left\|N^k\left(\eta_1\right)-N^k\left(\eta_2\right)\right\|_{\mathcal{L}(H_*^{-1 / 2}, H_*^{1 / 2})} \leq k \cdot\left\|N^k\right\| \cdot\left\|\eta_1-\eta_2\right\|_{W^{1, \infty}} \cdot M^{k-1}.
		$$
		where $M=\max \left(\left\|\eta_1\right\|_{W^{1, \infty}},\left\|\eta_2\right\|_{W^{1, \infty}}\right)$.
		Using the bound $\left\|N^k\right\| \leq C B^k$ from the analyticity result, we get
		$$
		\left\|N^k\left(\eta_1\right)-N^k\left(\eta_2\right)\right\|_{\mathcal{L}(H_*^{-1 / 2}, H_*^{1 / 2})} \leq k C B^k M^{k-1}\left\|\eta_1-\eta_2\right\|_{W^{1, \infty}},
		$$
		Since the continuity of the operators $\partial_{x_i}: H^{1/2}(\mathbb R^2)\to H^{-1/2}_{*}(\mathbb R^2)$ and $\partial_{x_i}: H^{1/2}_{*}(\mathbb R^2) \to H^{-1/2}(\mathbb R^2)$, we have
		$$
		\left\|L^k\left(\eta_1\right)-L^k\left(\eta_2\right)\right\|_{\mathcal{L}(H^{1/2}, H^{-1/2})} \leq k C B^k M^{k-1}\left\|\eta_1-\eta_2\right\|_{W^{1, \infty}}.
		$$
		Summing over $k$ :
		$$
		\|L(\eta_1)-L(\eta_2)\|_{\mathcal{L}(H^{1/2}, H^{-1/2})} \leq C\|\eta_1-\eta_2\|_{W^{1, \infty}} \sum_{k=1}^{\infty} k B^k M^{k-1}.
		$$
		If $M<\frac{1}{2 B}$, the series converges, and there exists $C>0$ such that:
		$$
		\left\|L\left(\eta_1\right)-L\left(\eta_2\right)\right\|_{\mathcal{L}(H^{1/2}, H^{-1/2})} \leq C\left\|\eta_1-\eta_2\right\|_{W^{1, \infty}}.
		$$
		This shows that $L(\cdot)$ is Lipschitz continuous in $B_{R}(0)$ with the $W^{1, \infty}$-norm.  To obtain Hölder continuity in the $H^1$-norm, we assume $\eta \in H^3\left(\mathbb{R}^2\right)$ (to control the $W^{1, \infty}$-norm). Using the Sobolev interpolation inequality (between $H^1$ and $H^3$ ):
		$$
		\|\eta\|_{W^{1, \infty}} 
		\leq C_{\delta}\|\eta\|_{H^1}^{1-\delta}\|\eta\|_{H^3}^{\delta}, \quad \delta \in(0,1) .
		$$
		Assume $\left\|\eta_1\right\|_{H^3},\left\|\eta_2\right\|_{H^3} \leq R$. Then
		$$
		\left\|\eta_1-\eta_2\right\|_{W^{1, \infty}} \leq C_{\delta}\left\|\eta_1-\eta_2\right\|_{H^1}^{1-\delta}(2 R)^{\delta}.
		$$
		We obtain
		$$
		\left\|L\left(\eta_1\right)-L\left(\eta_2\right)\right\|_{\mathcal{L}(H^{1/2}, H^{-1/2})} \leq  C_{\delta}(2 R)^{\delta}\left\|\eta_1-\eta_2\right\|_{H^1}^{1-\delta}
		$$
		Let $\alpha=1-\delta \in(0,1)$. Then, we have
		$$
		\left\|L\left(\eta_1\right)-L\left(\eta_2\right)\right\|_{\mathcal{L}(H^{1/2}, H^{-1/2})} \leq C\left\|\eta_1-\eta_2\right\|_{H^1}^\alpha
		$$
		The proof is thus complete.
	\end{proof}

	\begin{prop}\label{lemA5}
		There exist constants $\alpha \in (0,1)$ and $C > 0$ such that for all $\eta_1, \eta_2 \in \mathscr{E}_R$, 
		\begin{equation}\label{A5}
			\|K(\eta_1) - K(\eta_2)\|_{\mathcal{L}(H^{1/2}, H^{-1/2})} \le C \|\eta_1 - \eta_2\|^{\alpha}_{H^1(\mathbb{R}^2)}.
		\end{equation}
	\end{prop}
	
	\begin{proof}
		Following an argument similar to that in the proof of Proposition~\ref{lemA4}, we obtain the desired result.
	\end{proof}
	
	\begin{corollary}\label{propA3}
		There exists a constant $C > 0$ such that for all $\eta \in \mathscr{E}_R$, we have
		\begin{align*}
			\|K(\eta)\|_{\mathcal{L}(H^{1/2}(\mathbb{R}^2), H^{-1/2}(\mathbb{R}^2))} \le C, \qquad
			\|L(\eta)\|_{\mathcal{L}(H^{1/2}(\mathbb{R}^2), H^{-1/2}(\mathbb{R}^2))} \le C.
		\end{align*}
	\end{corollary}
	
	\begin{proof}
		Applying Proposition~\ref{lemA4} and Proposition~\ref{lemA5} with $\eta_1 = \eta$ and $\eta_2 = 0$ yields the desired bounds.
	\end{proof}

	With the analysis of the locked inertia tensor in hand, we now turn to the three-dimensional water-wave problem via the Hamiltonian system
	\begin{equation}\label{Hamilton}
		\begin{cases}
			\eta_{t}=D_{\xi} \mathcal{H}_{c}(\eta,\xi)\\
			\xi_{t} =-D_{\eta} \mathcal{H}_{c}(\eta,\xi).
		\end{cases}
	\end{equation}
	
	We fix a particular wave speed $c \in (\frac{1}{\sqrt{1+\varepsilon_0}}, 1)$ and investigate the stability of the solitary wave $\bar{u}(c) = (\bar{\eta}(c), \bar{\xi}(c))$. Linearizing \eqref{Hamilton} around $\bar{u}(c)$ with $X = (v,w)^T = u - \bar{u}(c)$ yields
	\begin{equation}\label{5-6}
		X_{t} = J\mathcal{A} X, \quad \text{where } J = \begin{pmatrix}
			0 & I \\
			-I & 0
		\end{pmatrix}, \quad \mathcal{A} = D^2\mathcal{H}_{c}(\bar{u}(c)).
	\end{equation}

	Using \eqref{2.3}, the self-adjoint operator $\mathcal{A}$ can be expressed as
	\begin{equation}\label{A-decom}
		\langle \langle \mathcal{A}X, X\rangle \rangle = \langle \bar{C} v, v\rangle + \langle \bar{G} (w + \bar{B} v), w + \bar{B} v \rangle,
	\end{equation}
	where $\langle \langle X_1, X \rangle \rangle := \langle v_1, v\rangle + \langle w_1, w\rangle$ for $X_1 = (v_1,w_1)^T \in \mathscr{F}^* = H^{-1}(\mathbb{R}^2) \times H_{*}^{-1/2}(\mathbb{R}^2)$ and $X=(v,w)^T \in \mathscr{F}=H^1(\mathbb R^2)\times H^{1/2}_{*}(\mathbb R^2)$, with 
	\begin{align*}
		\bar{C} &= D^2\mathcal{V}_{c}^{aug}(\bar{\eta}(c)), \qquad 
		\bar{G} = G(\bar{\eta}(c)), \\
		\bar{B} v &= \bar{G}^{-1} \mathcal{L}_{\bar{u}(c)} v = c \bar{G}^{-1} \left[\partial_{x_1}v - DG(\bar{\eta})[v] \bar{G}^{-1} \partial_{x_1} \bar{\eta}\right] \in \mathcal{L}(H^1(\mathbb{R}^2), H_{*}^{1/2}(\mathbb{R}^2)).
	\end{align*}

	Since $\bar{u}(c)$ is a traveling wave, for any $a, b \in \mathbb{R}$, the function $\bar{u}_{a,b}(c)(x_1,x_2) = \bar{u}(c)(x_1+a, x_2+b)$ is also a solution, and $\mathcal{H}_{c}(\bar{u}_{a,b}(c)) = \mathcal{H}_{c}(\bar{u}(c))$. From $D\mathcal{H}_{c}(\bar{u}(c)) = 0$, it follows that $D\mathcal{H}_{c}(\bar{u}_{a,b}(c)) = 0$. Differentiating with respect to $a$ at $a=0, b=0$ gives
	\begin{align*}
		\frac{d}{d a}D\mathcal{H}_{c}(\bar{u}_{a,b}(c)) \big|_{a=0,b=0} = D^2 \mathcal{H}_c(\bar{u}(c)) \partial_{x_1}\bar{u}(c) = \mathcal{A} \partial_{x_1}\bar{u}(c) = 0.
	\end{align*}
	This shows that $\partial_{x_1} \bar{u}(c)$ lies in the kernel of $D^2 \mathcal{H}_c(\bar{u}(c))$. Similarly, $\partial_{x_2} \bar{u}(c)$ is also in the kernel.

	Differentiating the identity $D\mathcal{H}_{c}(\bar{u}(c)) = 0$ with respect to the wave speed $c$ yields
	\[
	D^2\mathcal{H}_{c}(\bar{u}(c)) \partial_{c}\bar{u}(c) + \frac{\partial}{\partial c} D\mathcal{H}_{c}(\bar{u}(c)) = 0.
	\]
	Since $\frac{\partial}{\partial c} D\mathcal{H}_{c}(\bar{u}(c)) = D\mathcal{P}(\bar{u}(c)) = (-\partial_{x_1}\bar{\xi}(c), \partial_{x_1}\bar{\eta}(c))^T =: p_3$, we obtain
	\begin{equation}\label{5-7}
		\mathcal{A}(\partial_{x_1}\bar{u}(c)) = 0, \quad \mathcal{A}(\partial_{x_2}\bar{u}(c)) = 0, \quad \mathcal{A}(\partial_{c}\bar{u}(c)) = -p_3 = J(\partial_{x_1}\bar{u}(c)).
	\end{equation}
	This demonstrates that the vectors
	\begin{align*}
		z_1 := (\partial_{x_1}\bar{\eta}(c), \partial_{x_1}\bar{\xi}(c))^T, \quad
		z_2 := (\partial_{x_2}\bar{\eta}(c), \partial_{x_2}\bar{\xi}(c))^T, \quad
		z_3 &:= \partial_{c}\bar{u}(c)
	\end{align*}
	lie in the generalized kernel of $L = J\mathcal{A}$,  since $L z_1 = 0$, $L z_2 = 0$, and $L^2 z_3 = -L z_1=0$.

	Define $p_1 := (\partial_{x_1}\bar{\eta}(c), 0)^T$ and $p_2 := (\partial_{x_2}\bar{\eta}(c), 0)^T$ in $\mathscr{F}^*$. Let $\mathscr{F}^0 := \operatorname{span} \{z_1, z_2, z_3\}$ and
	\begin{equation*}
		\mathscr{F}^+ := \left\{ z = (v,w)^T \in \mathscr{F} \mid \langle\langle p_1, z \rangle\rangle = \langle\langle p_2, z \rangle\rangle = \langle\langle p_3, z \rangle\rangle = 0 \right\}.
	\end{equation*}

	\begin{prop}\label{lem5-3}
		The space $\mathscr{F}$ admits the direct sum decomposition $\mathscr{F} = \mathscr{F}^0 \oplus \mathscr{F}^+$, and the operator $\mathcal{A}$ exhibits the block structure
		\begin{equation}\label{5-8}
			\langle\langle \mathcal{A}(\textstyle{\sum_{i=1}^3 \alpha_i z_i + z_+}), (\textstyle{\sum_{i=1}^3 \alpha_i z_i + z_+}) \rangle\rangle = -d''(c) \alpha_3^2 + \langle\langle \mathcal{A}z_+, z_+ \rangle\rangle,
		\end{equation}
		for all $\alpha_i \in \mathbb{R}$, $i = 1,2,3$ and $z_+ \in \mathscr{F}^+$. Moreover, there exists a constant $\alpha_{+} > 0$ such that for all $z_+ \in \mathscr{F}^+$,
		\[
		\langle\langle \mathcal{A}z_+, z_+ \rangle\rangle \ge \alpha_+ \|z_+\|_\mathscr{F}^2.
		\]
	\end{prop}
	
	\begin{proof}
		We prove this proposition in three steps: (i) establishing the orthogonal decomposition of $\mathscr{F}$ and the block-diagonal form of $\mathcal{A}$, (ii) showing that $\mathcal{A}$ is positive definite up to a rank-three correction, and (iii) proving the coercivity of $\mathcal{A}$ on $\mathscr{F}^+$ via a contradiction argument.
		
		(i) First, we verify the orthogonality relations and construct the projection onto $\mathscr{F}^0$. Recall that $z_1 = (\partial_{x_1} \bar{\eta}(c), \partial_{x_1} \bar{\xi}(c))^T$, $z_2 = (\partial_{x_2} \bar{\eta}(c), \partial_{x_2} \bar{\xi}(c))^T$, and $z_3 = \partial_c \bar{u}(c)$. The dual elements are $p_1 = (\partial_{x_1} \bar{\eta}(c), 0)^T$, $p_2 = (\partial_{x_2} \bar{\eta}(c), 0)^T$, and $p_3 =  (-\partial_{x_1} \bar{\xi}(c), \partial_{x_1} \bar{\eta}(c))^T$.
		Using the relations $c = (1 + \ep^2)^{-1/2}$ and $\bar{\eta}_\ep = \ep^2 h(\ep x_1, \ep^2 x_2)$ (from the small-amplitude expansion in Section \ref{sect3}), we compute
		$$\partial_c \bar{\eta}(c) = -(1 + \ep^2)^{3/2} \ep^{-1} \partial_\ep \bar{\eta}_\ep = -(1 + \ep^2)^{3/2} (2 h + \ep x_1 \cdot \partial_1 h + 2 \ep^2 x_2 \cdot \partial_2 h).$$
		Given the symmetries of the solitary wave $\bar{\eta}(c)$ is even in both $x_1$ and $x_2$, while $\bar{\xi}(c)$ is odd in $x_1$ and even in $x_2$, it follows that
		$$\langle\langle p_i, z_j \rangle\rangle = 0 \quad \text{for } i \neq j.$$
		Moreover,
		$$\langle\langle p_i, z_i \rangle\rangle = \int_{\mathbb{R}^2} (\partial_{x_i} \bar{\eta})^2 \, dx > 0 \quad \text{for } i = 1, 2,$$
		since $\partial_{x_i} \bar{\eta} \not\equiv 0$. For the third element,
		\[
		\langle\langle p_3, z_3 \rangle\rangle = \langle\langle -\mathcal{A} z_3, z_3 \rangle\rangle = \langle\langle D\mathcal{P}(\bar{u}(c)), \partial_{c}\bar{u}(c) \rangle\rangle = d''(c) > 0.
		\]
		where the positivity of $d''(c)$ follows from the properties of the solitary waves as established in Proposition \ref{prop3.2}.
		
		Now define the projection operator $P: \mathscr{F} \to \mathscr{F}^0$ by
		\[
		Pz = \sum_{i=1}^3 \frac{\langle\langle p_i, z \rangle\rangle}{\langle\langle p_i, z_i \rangle\rangle} z_i.
		\]
		This projection has range $\mathscr{F}^0$ and kernel $\mathscr{F}^+ = \{ z \in \mathscr{F} \mid \langle\langle p_i, z \rangle\rangle = 0 \text{ for } i=1,2,3 \}$, establishing the direct sum decomposition $\mathscr{F} = \mathscr{F}^0 \oplus \mathscr{F}^+$. Using the relations in \eqref{5-7} and the self-adjointness of $\mathcal{A}$, we compute
		$$\langle\langle \mathcal{A} z, z \rangle\rangle = -d''(c) \alpha_3^2 + \langle\langle \mathcal{A} z_+, z_+ \rangle\rangle$$
		for $z = \sum_{i=1}^3 \alpha_i z_i + z_+$, confirming the block structure.
		
		(ii) Next, we show that $\mathcal{A}$ is positive definite up to a rank-three correction. Recall that $\bar{C} = D^2 \mathcal{V}_c^{{aug}}(\bar{\eta}(c))$. From the spectral analysis in Proposition \ref{prop3.2},  the operator $\bar{C}$ has a simple negative eigenvalue $\sigma_1(c)$ corresponding to an eigenfunction $\psi_1$, the eigenvalue $\sigma_2(c)=\sigma_3(c)=0$ associated with $\partial_{x_1}\bar{\eta}(c), \partial_{x_2}\bar{\eta}(c)$, and is strictly positive definite on the $L^2$-orthogonal complement of these three eigenspaces. Thus, for any $\delta > 0$, the perturbed operator
		\[
		\bar{C} + (\delta - \sigma_{1}(c)) \psi_1 \otimes \psi_1 + \delta (\partial_{x_1}\bar{\eta}) \otimes (\partial_{x_1}\bar{\eta}) + \delta (\partial_{x_2}\bar{\eta}) \otimes (\partial_{x_2}\bar{\eta})
		\]
		is positive definite on $L^2(\mathbb R^2)$ (the rank-three correction accounts for the negative direction and the translational invariances).  
		
		Since $\bar{G} = G(\bar{\eta}(c))$ is positive definite on $H^{1/2}_*(\mathbb{R}^2)$ by Lemma \ref{DNO bound}, combining with the decomposition \eqref{A-decom},
		$$\langle\langle \mathcal{A} z, z \rangle\rangle = \langle \bar{C} v, v \rangle + \langle \bar{G} (w + \bar{B} v), w + \bar{B} v \rangle,\quad z=(v,w)^{T},$$
		there exists $\alpha_\delta > 0$ such that
		\begin{align*}
			\langle\langle \mathcal{A}z, z \rangle\rangle &+ (\delta - \sigma_1(c)) \langle \langle (\psi_1, 0)^T, z \rangle \rangle^2 \\
			&+ \delta \langle \langle (\partial_{x_1}\bar{\eta}, 0)^T, z \rangle \rangle^2 + \delta \langle \langle (\partial_{x_2}\bar{\eta}, 0)^T, z \rangle \rangle^2 \ge \alpha_{\delta} \|z\|_{\mathscr{F}}^2
		\end{align*}
		for all $z \in \mathscr{F}$.

		(iii) Finally,  we establish the coercivity bound on $\mathscr{F}^+$. If not, there exist $z_4 \in \mathscr{F}^{+}$ such that $\left\langle\left\langle \mathcal{A} z_4, z_4\right\rangle\right\rangle \leq 0$. By (i), we know that $\mathcal{A}$ is negative semi-definite on $\mathscr{F}^0$. Then it follows that $\mathcal{A}$ is negative semi-definite on $\left\{z_4\right\} \oplus \mathscr{F}^0$, which is contradicted by (ii).
	\end{proof}

	\begin{theorem}
		Every solitary capillary-gravity water wave $\bar{u}(c)=(\bar{\eta}(c), \bar{\xi}(c))$ defined in \eqref{bareta} with wave speed $c = \frac{1}{\sqrt{1+\varepsilon^2}}$ for $\varepsilon \in (0,\varepsilon_0)$ is conditionally orbitally stable in the following sense: For every $R > 1$ and $\rho > 0$, there exists $\rho_0 > 0$ such that if 
		$u=(\eta,\xi) \colon [0,T) \to \mathscr{F}_R$ is a continuous solution of \eqref{Hamilton} that preserves the functionals $\mathcal{H}$ and $\mathcal{P}$, and if the initial data satisfy
		\[
		\|\eta(0) - \bar{\eta}(c)\|_{H^1(\mathbb{R}^2)} + \|\xi(0) - \bar{\xi}(c)\|_{H^{1/2}_{*}(\mathbb{R}^2)} < \rho_0,
		\]
		then for all $t \in [0,T)$,
		\begin{equation*}
			\inf_{(x_0,y_0) \in \mathbb{R}^2} \left( \|\eta(t, \cdot - (x_0,y_0)) - \bar{\eta}(c)\|_{H^1(\mathbb{R}^2)} + \|\xi(t, \cdot - (x_0,y_0)) - \bar{\xi}(c)\|_{H^{1/2}_{*}(\mathbb{R}^2)} \right) < \rho.
		\end{equation*}
	\end{theorem}

	\begin{proof}
		Recall that \(\mathcal{P}\) is a conserved quantity for the water wave system (due to the Hamiltonian structure in \eqref{1-1}), so \(\mathcal{P}(u(t)) = \mathcal{P}(u(0))\) for all \(t \in \mathbb{R}\).
		
		Let \(c_0 = \frac{1}{\sqrt{1+\varepsilon_0^2}} < 1\) be the initial wave speed such that the initial data \(u(0)\) is a small perturbation of \(\bar{u}(c_0)\) in the space \(\mathscr{F}\), with \(\|u(0) - \bar{u}(c_0)\|_\mathscr{F} \leq \rho_0\) for sufficiently small \(\rho_0 > 0\).
		
		Define the mapping \(p\colon c \mapsto \mathcal{P}(\bar{u}(c))\). For the solitary wave family \(\bar{u}(c) = (\bar{\eta}(c), \bar{\xi}(c))\) parameterized by \(c = (1+\varepsilon^2)^{-1/2}\), we compute the derivative:
		\[
		\partial_c p(c) = \partial_c \mathcal{P}(\bar{u}(c)) = D\mathcal{P}(\bar{u}(c))[\partial_c \bar{u}(c)] = d''(c).
		\]
		From Proposition~\ref{prop3.2}, we have the strict convexity \(d''(c) > 0\), and more precisely, the asymptotic estimate \(d''(c) > C \) holds for some constant \(C > 0\).
		
		Since 
		$$\|\eta(0)\|_{H^1(\mathbb{R}^2)} \le \rho_0 + \|\bar{\eta}(c)\|_{H^1(\mathbb{R}^2)} \le C$$ 
		and 
		$$\|\xi(0)\|_{H^{1/2}_{*}(\mathbb{R}^2)} \le \rho_0 + \|\bar{\xi}(c)\|_{H^{1/2}_{*}(\mathbb{R}^2)} \le C,$$ 
		we obtain
		\[
		|\mathcal{P}(u(0)) - \mathcal{P}(\bar{u}(c_0))| \leq C\|u(0) - \bar{u}(c_0)\|_\mathscr{F} \leq C\rho_0.
		\]
		Since \(p(c)\) is strictly increasing and continuous, the Inverse Function Theorem guarantees the existence of a unique \(c\) in a neighborhood of \(c_0\) such that
		\[
		\mathcal{P}(\bar{u}(c)) = \mathcal{P}(u(0)).
		\]
		Moreover, applying the mean value theorem yields the estimate
		\[
		|c - c_0| \leq \frac{|\mathcal{P}(u(0)) - \mathcal{P}(\bar{u}(c_0))|}{\inf_{c \in I} p'(c)} \leq \frac{C\rho_0}{\inf_{c \in I} d''(c)} \leq C\rho_0,
		\]
		where \(I\) is a sufficiently small interval containing \(c_0\).
		Thus, by redefining \(c \in (c_0 - C\rho_0, c_0 + C\rho_0)\) for sufficiently small \(\rho_0\), we ensure that
		\[
		\mathcal{P}(u(t)) = \mathcal{P}(u(0)) = \mathcal{P}(\bar{u}(c)) \quad \text{for all } t \in \mathbb{R}.
		\]
		We begin with the expression for the augmented Hamiltonian at the initial data:
		\begin{align*}
			\mathcal{H}_{c}(u(0)) &= \mathcal{V}(\eta(0)) + \mathcal{K}(\eta(0),\xi(0)) + c\mathcal{P}(\eta(0),\xi(0)).
		\end{align*}
		Since $\bar{u}(c) = (\bar{\eta}(c), \bar{\xi}(c))$ is a critical point of $\mathcal{H}_c$, we have $D\mathcal{H}_c(\bar{u}(c)) = 0$. We consider the Taylor expansion of $\mathcal{H}_c$ around $\bar{u}(c)$
		\begin{align*}
			\mathcal{H}_c(u(0)) &= \mathcal{H}_c(\bar{u}(c)) + \frac{1}{2} D^2\mathcal{H}_c(\bar{u}(c))[u(0) - \bar{u}(c), u(0) - \bar{u}(c)] \\
			&\quad + R(u(0) - \bar{u}(c)),
		\end{align*}
		where $R(u(0) - \bar{u}(c))$ denotes the remainder term of order $O(\|u(0) - \bar{u}(c)\|_\mathscr{F}^3)$. 
		From the expression for the second variation in equation \eqref{2.3}, we have
		\begin{align*}
			D^2\mathcal{H}_c(\bar{u}(c))\left[\begin{pmatrix}v\\ w\end{pmatrix},\begin{pmatrix}v\\ w\end{pmatrix}\right] 
			&= D^2\mathcal{V}(\bar{\eta}(c))[v,v] + \frac{1}{2}\langle D^2G(\bar{\eta}(c))[v,v]\bar{\xi}(c), \bar{\xi}(c)\rangle \\
			&\quad + \langle G(\bar{\eta}(c))w, w\rangle + 2\langle \mathcal{L}_{\bar{u}(c)}v, w\rangle,
		\end{align*}
		where $v = \eta(0) - \bar{\eta}(c)$, $w = \xi(0) - \bar{\xi}(c)$, and $\mathcal{L}_{\bar{u}(c)}v = c\partial_{x_1}v + DG(\bar{\eta}(c))[v]\bar{\xi}(c)$.
		To estimate this quadratic form, we use the following bounds
		\begin{enumerate}
			\item By Lemma \ref{DNO bound}, the Dirichlet-Neumann operator satisfies
			\[
			c_1\|w\|_{*,1/2}^2 \leq \langle G(\bar{\eta}(c))w, w\rangle \leq c_2\|w\|_{*,1/2}^2.
			\]
			\item The potential energy term is bounded by
			\[
			|D^2\mathcal{V}(\bar{\eta}(c))[v,v]| \leq C\|v\|_{H^1}^2.
			\]
			\item By Lemma \ref{DG},  the term involving $D^2G(\bar{\eta}(c))$ satisfies
			\[
			\langle D^2G(\bar{\eta}(c))[v,v]\bar{\xi}(c), \bar{\xi}(c)\rangle = \langle v, a_3 v \rangle+2\langle a_2v, G(\bar{\eta}(c))a_2v\rangle.
			\]
			Since $\bar{\eta}(c)$ and $\bar{\xi}(c)$ are fixed solitary wave solutions with sufficient regularity ($\bar{\eta}(c) \in H^s(\mathbb{R}^2)$ with $s > 5/2$, $\bar{\xi}(c) \in H^{1/2}_*(\mathbb{R}^2)$), the coefficients $a_1$ and $a_2$ are bounded in $W^{1,\infty}(\mathbb{R}^2)$ with bounds depending on $\|\bar{\eta}(c)\|_{W^{2,\infty}}$ and $\|\bar{\xi}(c)\|_{*,1/2}$.
			Therefore, we have the estimates
			\[
			|\langle v, a_3 v \rangle| \leq \|a_3\|_{L^\infty} \|v\|_{L^2}^2 \leq C \|v\|_{H^1}^2,
			\]    
			and by Lemma \ref{DNO bound}
			\[
			|\langle a_2 v, G(\bar{\eta}(c))(a_2 v) \rangle| \leq c_2 \|a_2 v\|_{*,1/2}^2 \leq C \|v\|_{H^1}^2.
			\]
			
			\item For the cross term, we have
			\[
			|2\langle \mathcal{L}_{\bar{u}(c)}v, w\rangle| \leq  \|w\|_{*,1/2}^2 + \|\mathcal{L}_{\bar{u}(c)}v\|_{*,-1/2}^2.
			\]
			Now estimate $\|\mathcal{L}_{\bar{u}(c)}v\|_{*,-1/2}$:
			\[
			\|\mathcal{L}_{\bar{u}(c)}v\|_{*,-1/2} \leq |c| \|\partial_{x_1}v\|_{*,-1/2} + \|DG(\bar{\eta}(c))[v]\bar{\xi}(c)\|_{*,-1/2}.
			\]
			For the first term, since $\partial_{x_1}: H^{1/2}(\mathbb{R}^2) \to H^{-1/2}_*(\mathbb{R}^2)$ is bounded
			\[
			\|\partial_{x_1}v\|_{*,-1/2} \leq C \|v\|_{1/2} \leq C \|v\|_{H^1}.
			\]
			For the second term, by Lemma \ref{DG} and the definition of the dual norm
			\[
			\|DG(\bar{\eta}(c))[v]\bar{\xi}(c)\|_{*,-1/2} = \sup_{\|\psi\|_{*,1/2}=1} |\langle \psi, DG(\bar{\eta}(c))[v]\bar{\xi}(c) \rangle|.
			\]
			Using the representation from Lemma \ref{DG}
			\[
			|\langle \psi, DG(\bar{\eta}(c))[v]\bar{\xi}(c) \rangle| = |\langle v, a_1 \cdot \nabla \psi + a_2 G(\bar{\eta}(c))\psi \rangle|.
			\]
			Since $a_1, a_2 \in W^{1,\infty}(\mathbb{R}^2)$, $\partial_{x_i}: H^{1/2}(\mathbb{R}^2) \to H^{-1/2}_*(\mathbb{R}^2)$ is bounded for $i=1,2$,  and $G(\bar{\eta}(c))$ is bounded from $H^{1/2}_*$ to $H^{-1/2}_*$
			\begin{align*}
				|\langle v, a_1 \cdot \nabla \psi \rangle|& =|\langle \div(a_1v),   \psi \rangle|\le \|\div(a_1v)\|_{*,-1/2}  \|\psi\|_{*,1/2}\\
				&\le C\|a_1v\|_{H^{1/2}}\|\psi\|_{*,1/2}\le C\|v\|_{H^1}\|\psi\|_{*,1/2},
			\end{align*}
			and
			\begin{align*}
				|\langle v, a_2 G(\bar{\eta}(c))\psi \rangle| & \leq \sqrt{\langle a_2 v, G(\bar{\eta}(c)) (a_2 v)\rangle } \sqrt{\langle \psi, G(\bar{\eta}(c)) \psi\rangle}\\
				&\le C\|a_2v\|_{*, 1/2}\|\psi\|_{*, 1/2}\le C\|v\|_{H^1}\|\psi\|_{*, 1/2}.
			\end{align*}
			Therefore, we have
			\[
			\|DG(\bar{\eta}(c))[v]\bar{\xi}(c)\|_{*,-1/2} \leq C \|v\|_{H^1},
			\]
			and consequently,
			\[
			\|\mathcal{L}_{\bar{u}(c)}v\|_{*,-1/2}^2 \leq C\|v\|_{H^1}^2.
			\]
		\end{enumerate}
		Combining these estimates, we obtain
		\[
		|D^2\mathcal{H}_c(\bar{u}(c))[u(0) - \bar{u}(c), u(0) - \bar{u}(c)]| \leq K_1\|u(0) - \bar{u}(c)\|_\mathscr{F}^2,
		\]
		for some constant $K_1 > 0$.
		For the remainder term, by the analyticity of $\mathcal{H}_c$ (established through the analyticity of $G(\eta)$ and $G^{-1}(\eta)$ in Lemmas \ref{DNO bound} and \ref{G inverse analytic}), we have
		\[
		|R(u(0) - \bar{u}(c))| \leq C\|u(0) - \bar{u}(c)\|_\mathscr{F}^3.
		\]
		Therefore, for sufficiently small $\|u(0) - \bar{u}(c)\|_\mathscr{F} \leq \rho_0$, we conclude
		\begin{align*}
			\mathcal{H}_{c}(u(0)) &\leq \mathcal{H}_{c}(\bar{u}(c)) + K_1\|u(0) - \bar{u}(c)\|_\mathscr{F}^2 + C\rho_0^3 \\
			&\leq \mathcal{H}_{c}(\bar{u}(c)) + (K_1 + C\rho_0)\rho_0^2 \\
			&\leq \mathcal{H}_{c}(\bar{u}(c)) + K\rho_0^2,
		\end{align*}
		where $K = K_1 + C\rho_0$ is a positive constant.
		Since the flow preserves both $\mathcal{H}$ and $\mathcal{P}$, it preserves $\mathcal{H}_c = \mathcal{H} + c \mathcal{P}$. Thus:
		$$\mathcal{H}_c(u(t)) = \mathcal{H}_c(u(0)) \leq \mathcal{H}_c(\bar{u}(c)) + K\rho_0^2, \quad \forall t \in [0,T).
		$$

		We now proceed to construct the optimal translation and establish orbital stability via a continuity argument.

		\paragraph{\bf Step 1: Local existence of optimal translation}
		Define $F: \mathbb R^2 \times \mathscr{F} \mapsto \mathbb R^2$, $F=(F_1, F_2)^T$ with 
		\[
		F_i(a,b, u)=\int_{\mathbb R^2} (\partial_{x_i}\bar \eta ) (\eta(x_1-a, x_2-b)-\bar \eta )dx_1dx_2,
		\]
		for $i=1,2$ and $u(t)=(\eta(t), \xi(t))$. We then obtain that $F$ is continuously differentiable and $F(0,0,\bar u)=(0,0)$. By direct calculation, we have 
		\begin{align*}
			\frac{\partial F_1}{\partial a}\Big|_{(0,0,\bar u)}=-\int_{\mathbb R^2} (\partial_{x_1}\bar \eta)^2dx_1dx_2<0,\quad \frac{\partial F_2}{\partial b}\Big|_{(0,0,\bar u)}=-\int_{\mathbb R^2} (\partial_{x_2}\bar \eta)^2dx_1dx_2<0,
		\end{align*}
		and by symmetry of $\bar \eta$,
		\begin{align*}
			\frac{\partial F_1}{\partial b}\Big|_{(0,0,\bar u)}=0,\quad \frac{\partial F_2}{\partial a}\Big|_{(0,0,\bar u)}=0.
		\end{align*}
		Then the Jacobi matrix $D_{(a,b)}F(0,0,\bar u)$ is  invertible.  
		By the implicit function theorem, there exists $\delta > 0$ such that for every $u \in B_{\delta}(\bar u)$, there is a unique smooth function $u \mapsto (a_{}(u), b_(u)) \in \mathbb R^2$ satisfying
		\[
		F(a_*(u), b_*(u), u)=(0,0)\quad  \text{with} \quad  (a_*(\bar u), b_*(\bar u))=(0,0).
		\]

		\paragraph{\bf Step 2: Continuity argument for global existence}
		Define the orbital distance
		$$\mu(t) := \inf_{(x_0, y_0) \in \mathbb{R}^2} \| u(t, \cdot - (x_0, y_0)) - \bar{u}(c) \|_\mathscr{F}.$$
		Note that $\mu(t)$ is continuous in $t$ due to the continuity of $u(t)$ in $\mathscr{F}$.
		Choose $\sigma > 0$ such that $\sigma < \delta$ and define the maximal time
		$$T_* = \sup \left\{ t \in [0,T) : \mu(s) < \sigma \text{ for all } 0 \leq s \leq t \right\}.$$
		We claim that $T_* = T$. Suppose for contradiction that $T_* < T$. Then by continuity of $\mu(t)$, we have $\mu(T_*) = \sigma$.
		
		For all $t \in [0, T_*)$, since $\mu(t) < \sigma$, there exists a rough translation $(x_0(t), y_0(t))$ such that the roughly translated solution
		$$
		v(t) := u(t, \cdot - (x_0(t), y_0(t)))
		$$
		satisfies $\|v(t) - \bar{u}\|_\mathscr{F} < \sigma <  \delta$, so $v(t) \in B_{\delta}(\bar{u})$. Thus, the optimal fine translation $(a_*(v(t)), b_*(v(t)))$ is well-defined by step 1, and we define
		$$
		z(t) := v(t, \cdot - (a_*(v(t)), b_*(v(t)))) - \bar{u}.
		$$
		From the construction, we have $\langle \langle p_i, z(t) \rangle \rangle = 0$ for $i=1,2$, where $p_i = (\partial_{x_i}\bar{\eta}, 0)$.

		Moreover, translations preserve the conservation laws, so
		$$\mathcal{P}(\bar{u} + z(t)) = \mathcal{P}(\bar{u}), \quad \mathcal{H}_c(\bar{u} + z(t)) \leq \mathcal{H}_c(\bar{u}) + K\rho_0^2.$$
		By Proposition \ref{lem5-3}, we have the decomposition $z(t) = \alpha_1 z_1 + \alpha_2 z_2 + \alpha_3 z_3 + z_+$. Since $\langle \langle p_i, z(t) \rangle \rangle = 0$ for $i=1,2$, we get $\alpha_1 = \alpha_2 = 0$. From $\langle \langle p_3,z_3 \rangle \rangle= d''(c)>C>0$ and Taylor expansion
		\[
		\mathcal{P}(\bar{u} + z(t)) = \mathcal{P}(\bar{u}) + \langle \langle D\mathcal{P}(\bar{u}), z(t) \rangle \rangle + O(\|z(t)\|_\mathscr{F}^2),
		\]
		we obtain $\langle \langle p_3, z(t) \rangle \rangle = O(\|z(t)\|_\mathscr{F}^2)$ and thus $\alpha_3 = O(\|z(t)\|_\mathscr{F}^2)$.
		
		For $z(t)=(v,w)^T$ and $\bar B v= \bar G ^{-1} \mathcal{L}_{\bar u (c)} v=c \bar G ^{-1} [\partial_{x_1}v-DG(\bar \eta)[v] \bar G ^{-1} \partial_{x_1} \bar \eta]$, we define $\mathcal{B}(v)=c  G(\bar \eta+v) ^{-1} \partial_{x_1}(\bar \eta+v)- c  G(\bar \eta) ^{-1} \partial_{x_1}(\bar \eta)$ that $\bar B$ is the linearization of $\mathcal{B}$ at $v=0$.

		Recall that $\mathcal{H}_c(\eta,\xi)=\mathcal{V}(\eta)+\frac{1}{2}\langle \xi, G(\eta)\xi \rangle +c\langle \xi, \partial_{x_1}\eta \rangle$. For fixed $\eta$, this can be rewritten by completing the square in $\xi$:
		$$\mathcal{H}_c(\eta, \xi) = \mathcal{V}_c^{{aug}}(\eta) + \frac{1}{2} \left\langle \xi + c G(\eta)^{-1} \partial_{x_1}\eta, G(\eta) (\xi + c G(\eta)^{-1} \partial_{x_1}\eta) \right\rangle,$$
		where $\mathcal{V}_c^{{aug}}(\eta) = \mathcal{V}(\eta) - \frac{1}{2} c^2 \int_{\mathbb{R}^2} \partial_{x_1}\eta [G(\eta)^{-1} \partial_{x_1}\eta] \, dx_1dx_2$ is the augmented potential, and the quadratic term is nonnegative (zero at the minimizing $\xi = -c G(\eta)^{-1} \partial_{x_1}\eta$). Substituting the perturbed values $\eta = \bar{\eta} + v$ and $\xi = \bar{\xi} + w$ gives
		$$\mathcal H_c(\bar{\eta} + v, \bar{\xi} + w) = \mathcal V_c^{\mathrm{aug}}(\bar{\eta} + v) + \frac{1}{2} \left\langle w + \mathcal B(v), G(\bar{\eta} + v) (w + \mathcal B(v)) \right\rangle,$$
		where $\mathcal B(v) = c G(\bar{\eta} + v)^{-1} \partial_{x_1}(\bar{\eta} + v) - c G(\bar{\eta})^{-1} \partial_{x_1}\bar{\eta}$ (noting $\bar{\xi} = -c G(\bar{\eta})^{-1} \partial_{x_1}\bar{\eta}$, so $\mathcal{B}(v) = \bar{\xi} + c G(\bar{\eta} + v)^{-1} \partial_{x_1}(\bar{\eta} + v)$).
		Then we derive 
		$$\mathcal H_c(\bar u + z) - \mathcal H_c(\bar u) = \left[ \mathcal V_c^{\mathrm{aug}}(\bar{\eta} + v) - \mathcal V_c^{\mathrm{aug}}(\bar{\eta}) \right] + \frac{1}{2} \left\langle w + \mathcal{B}(v), G(\bar{\eta} + v) (w + \mathcal{B}(v)) \right\rangle.$$
		The second variation operator $\mathcal{A} = D^2 \mathcal H_c(\bar{u})$ has quadratic form
		$$\frac{1}{2} \langle\langle \mathcal{A} z, z \rangle \rangle= \frac{1}{2} \langle \bar{C} v, v \rangle + \frac{1}{2} \left\langle w + \bar{B} v, \bar{G} (w + \bar{B} v) \right\rangle,$$
		where $\bar{C} = D^2 \mathcal V_c^{{aug}}(\bar{\eta})$, $\bar{G} = G(\bar{\eta})$, and $\bar{B} v = G(\bar{\eta})^{-1} \mathcal L_{\bar{u}(c)} v$ with $\mathcal L_{\bar{u}(c)} v = c \partial_{x_1}v + DG(\bar{\eta})[v] \bar{\xi}$. We have 
		\[
		K\rho_0^2\ge \mathcal{H}_c(\bar u+z)-\mathcal{H}_c(\bar u)=: J_1+J_2+J_3+J_4+J_5,
		\]
		where 
		\begin{align*}
			J_1 &= \tfrac{1}{2} \langle\langle \mathcal{A} z, z \rangle \rangle,\\
			J_2 &=  \mathcal V_c^{\mathrm{aug}}(\bar{\eta} + v) - \mathcal V_c^{\mathrm{aug}}(\bar{\eta})-\tfrac{1}{2} \langle \bar{C} v, v \rangle, \\
			J_3 &= \tfrac{1}{2} \langle G(\bar \eta+v)(\mathcal{B}(v)- \bar B v), w+\mathcal{B}(v)\rangle, \\
			J_4 &=\tfrac{1}{2} \langle (G(\bar \eta+v)-G(\bar \eta))(w+\bar B v), w+\bar B v\rangle,\\
			J_5 &=\tfrac{1}{2} \langle \mathcal{B}(v)-\bar B v, G(\bar \eta+v)(w+\bar B v) \rangle.
		\end{align*}
		For the term $J_1$, we derive from Proposition \ref{lem5-3} and  $\alpha_3=O(\|z\|_{\mathscr{F}}^2)$ that 
		\[
		J_1\ge \frac{1}{2} \alpha_{+}\|z_+\|_{\mathscr{F}}^2 -O(\|z\|_{\mathscr{F}}^4).
		\]
		For the second term $J_2$, we have $J_2=O(\|v\|^3_{H^1(\mathbb R^2)})$. For the term $J_3$, we notice from the positive definiteness of $ G(\bar \eta +v)$ on $H^{1/2}_{*}(\mathbb R^2)$ and Lemma \ref{DNO bound} that  the inner product $\langle G a, b \rangle$ satisfies the Cauchy-Schwarz inequality in the $ H^{1/2}_{*} $ norm:
		\begin{align*}
			&|\langle G(\bar \eta +v) (\mathcal{B}(v)- \bar B v), w+\mathcal{B}(v) \rangle|\\
			&\leq \sqrt{\langle G (\mathcal{B}(v)- \bar B v), (\mathcal{B}(v)- \bar B v) \rangle} \cdot \sqrt{\langle G (w+\mathcal{B}(v)), (w+\mathcal{B}(v)) \rangle}\\
			&\le  \|\mathcal{B}(v)- \bar B v\|_{*,1/2} \cdot \|w+\mathcal{B}(v)\|_{*,1/2}.
		\end{align*}
		We note that 
		\begin{align*}
			\partial_{x_1}(\mathcal{B}(v)- \bar B v)&=-c\left(K(\bar \eta+v)[\bar \eta+v]-K(\bar \eta)[\bar \eta]\right)-\partial_{x_1}\bar B v\\
			&=-c (K(\bar \eta+v)-K(\bar \eta))[v]\\
			&\quad -c\big( K(\bar \eta+v)[\bar \eta]-K(\bar \eta)[\bar \eta]- DK(\bar \eta)[v] \bar \eta  \big).
		\end{align*}
		Then, we obtain from Proposition \ref{lemA5} that  
		\begin{align*}
			\|\partial_{x_1}(\mathcal{B}(v)- \bar{B} v)\|_{H^{-1}(\mathbb R^2)}&\le C\left(\|K(\bar \eta+v)-K(\bar \eta)\|_{\mathcal{L}(H^1, H^{-1})}\|v\|_{H^1}+\|v\|_{H^1}^2\right)\\
			&\le C\left(\|K(\bar \eta+v)-K(\bar \eta)\|_{\mathcal{L}(H^{1/2}, H^{-1/2})}\|v\|_{H^1}+\|v\|_{H^1}^2\right) \\
			&\le C\|v\|_{H^1}^{1+\alpha}.
		\end{align*}
		Similarly, by Proposition \ref{lemA4}, we have
		\begin{align*}
			\|\partial_{x_2}(\mathcal{B}(v)- \bar B v)\|_{H^{-1}(\mathbb R^2)}&\le C\left(\|L(\bar \eta+v)-L(\bar \eta)\|_{\mathcal{L}(H^1, H^{-1})}\|v\|_{H^1}+\|v\|_{H^1}^2\right)\\
			&\le C\left(\|L(\bar \eta+v)-L(\bar \eta)\|_{\mathcal{L}(H^{1/2}, H^{-1/2})}\|v\|_{H^1}+\|v\|_{H^1}^2\right) \\
			&\le C\|v\|_{H^1}^{1+\alpha}.
		\end{align*}
		Thus, we have
		\begin{align*}
			\|\mathcal{B}(v)- \bar B v\|_{H^{1/2}_{*}(\mathbb R^2)}\le C\|\nabla(\mathcal{B}(v)- \bar B v)\|_{H^{-1}(\mathbb R^2)} \le C\|v\|_{H^1}^{1+\alpha}.
		\end{align*}
		For $\|w+\mathcal{B}(v)\|_{*,1/2}$, we obtain
		\[
		\|w+\mathcal{B}(v)\|_{*,1/2}\le \|w\|_{*,1/2}+\|\mathcal{B}(v)\|_{*,1/2}.
		\]
		Since $\bar B \in \mathcal{L}(H^1(\mathbb R^2), H_{*}^{1/2}(\mathbb R^2))$, we derive 
		\begin{align*} 
			\|\mathcal{B}(v)\|_{*,1/2}&\le \|\mathcal{B}(v)-\bar B v\|_{*,1/2}+\|\bar B v\|_{*,1/2}\le C\|v\|_{H^1}.
		\end{align*}
		Then we have
		\[
		|J_3|\le C\|v\|_{H^1}^{1+\alpha}\|z\|_{\mathscr{F}}.
		\]
		For the term $J_4$,  we derive from the analyticity of $G(\cdot)$ and similar argument in the proof of Proposition \ref{lemA4} that 
		\begin{align*}
			\left|J_4\right| & =\tfrac{1}{2}|\langle(G(\bar{\eta}+v)-\bar{G})(w+\bar{B} v), w+\bar{B} v\rangle| \\ 
			& \leq \tfrac{1}{2}\|G(\bar{\eta}+v)-\bar{G}\|_{H_*^{1 / 2} \rightarrow H_{*}^{-1 / 2}} \cdot\|w+\bar{B} v\|_{H_*^{1 / 2}}^2\\
			& \le C\|v\|_{H^1}^\alpha \|z\|_\mathscr{F}^2.
		\end{align*}
		For the term $J_5$, we obtain 
		\begin{align*}
			|J_5| &=\tfrac{1}{2} |\langle \mathcal{B}(v)-\bar B v, G(\bar \eta+v)(w+\bar B v) \rangle|\\
			&\le C  \sqrt{\langle  (\mathcal{B}(v)- \bar B v), G(\bar \eta+v)(\mathcal{B}(v)- \bar B v) \rangle} \cdot \sqrt{\langle (w+\bar{B} v), G(\bar \eta+v)(w+\bar{B} v) \rangle}\\
			&\le  C\|\mathcal{B}(v)- \bar B v\|_{*,1/2} \cdot \|w+\bar{B} v\|_{*,1/2}\\
			&\le C\|v\|_{H^1}^{1+\alpha}\|z\|_{\mathscr{F}}.
		\end{align*}
		Combining the above estimates for $J_1-J_5$, we have
		\begin{align*}
			K\rho_0^2\ge  \frac{1}{2} \alpha_{+}\|z_+\|_{\mathscr{F}}^2 -O(\|z\|_{\mathscr{F}}^4)+O(\|v\|^3_{H^1(\mathbb R^2)})+O(\|v\|_{H^1}^{1+\alpha}\|z\|_{\mathscr{F}})+O(\|v\|_{H^1}^\alpha \|z\|_\mathscr{F}^2).
		\end{align*}
		As $\|z_+\|_{\mathscr{F}}\le \|z\|_{\mathscr{F}}\le 2\|z_{+}\|_{\mathscr{F}}$ in a suitable neighbourhood of $z=0$, we obtain $K\rho_0^2\ge \alpha_{+}/9\|z\|_{\mathscr{F}}^2$ as long as $\|z\|_{\mathscr{F}}\le \rho_1$. 
		This implies that for all $t \in [0, T_*)$:
		\[
		\|z(t)\|_\mathscr{F} \leq \sqrt{\frac{9K}{\alpha_+}} \rho_0,
		\]
		provided we choose $\rho_0 < (\alpha_{+}/9K)^{1/2}\rho_1$.

		Now, since $\mu(t)$ is defined as the infimum over all translations, and $z(t)$ corresponds to a particular translation, we have
		\[
		\mu(t) \leq \|z(t)\|_\mathscr{F} \leq \sqrt{\frac{9K}{\alpha_+}} \rho_0.
		\]
		Choose $\rho_0 > 0$ sufficiently small such that
		\[
		\sqrt{\frac{9K}{\alpha_+}} \rho_0 < \min\{\rho, \sigma\}.
		\]
		Then, for $t < T_*$, $\mu(t) < \sigma$, and at $t = T_*$, the estimates would yield $\mu(T_*) < \sigma$, contradicting $\mu(T_*) = \sigma$. Therefore, $T_* = T$, and $\mu(t) < \sigma$ for all $t \in [0, T)$. 
		Finally, choosing $\rho_0 < (\alpha_{+}/9K)^{1/2} \min\{\rho, \rho_1, \sigma \}$ ensures that $\mu(t) < \rho$ for all $t \in [0, T)$.
	\end{proof}

	\appendix 
	\renewcommand{\thedefinition}{A.\arabic{definition}}
	\setcounter{definition}{0}
	\renewcommand{\thetheorem}{A.\arabic{theorem}}
	\setcounter{theorem}{0}
	\renewcommand{\thelemma}{A.\arabic{lemma}}
	\setcounter{lemma}{0}
	\renewcommand{\thecorollary}{A.\arabic{corollary}}
	\setcounter{corollary}{0}
	\begin{appendices}
		
		\section*{Appendix}

		In this appendix, we collect some fundamental definitions and technical results that are essential to the analysis in the main body of the paper. We begin by recalling the definition of analyticity for operators between Banach spaces, which plays a crucial role in our study of water wave problems. This is followed by several key lemmas concerning the analytic properties of the Dirichlet-Neumann operator and its inverse, which are central to our approach.

		\begin{definition}
			Let $X_1$ and $X_2$ be Banach spaces, $U$ be a non-empty, open subset of $X_1$, and $\mathcal{L}_{\mathrm{s}}^k(X_1, X_2)$ denote the space of bounded, symmetric $k$-linear operators from $X_1^k$ to $X_2$, equipped with the norm
			\[
			\|m\| := \inf \left\{ \gamma > 0 : \|m(\{f\}^{(k)})\|_{X_2} \le \gamma \|f\|_{X_1}^k \text{ for all } f \in X_1 \right\}.
			\]
			A function $F: U \to X_2$ is said to be analytic at a point $x_0 \in U$ if there exist constants $\delta, r > 0$ and a sequence of operators $\{m_k\}_{k=0}^\infty$, where $m_k \in \mathcal{L}_{\mathrm{s}}^k(X_1, X_2)$, such that
			\[
			F(x) = \sum_{k=0}^{\infty} m_k(\{x - x_0\}^{(k)}) \quad \text{for all } x \in B_\delta(x_0),
			\]
			and
			\[
			\sup_{k \ge 0} \, r^k \|m_k\| < \infty.
			\]
		\end{definition}

		This notion of analyticity provides a powerful framework for studying nonlinear operators through their power series expansions. In the context of water waves, it allows us to rigorously analyze the dependence of the Dirichlet-Neumann operator on the free surface elevation.
		
		Building upon this definition, we now present two fundamental results concerning the analytic behavior of boundary operators in water wave theory. The first lemma establishes the analyticity of the Dirichlet-Neumann operator and provides crucial bounds for its quadratic form.

		\begin{lemma}\label{DNO bound}
			Let $G(\cdot): W^{1,\infty}(\mathbb{R}^2) \to \mathcal{L}(H_*^{1/2}(\mathbb{R}^2), H_*^{-1/2}(\mathbb{R}^2))$ be the Dirichlet--Neumann operator. Then:
			\begin{enumerate}
				\item[(i)] $G(\cdot)$ is analytic at the origin.
				\item[(ii)] There exist constants $c_1, c_2 > 0$ and $R > 0$ such that for all $\eta \in B_R(0) \subset W^{1,\infty}(\mathbb{R}^2)$ and all $\Phi \in H_*^{1/2}(\mathbb{R}^2)$,
				\[
				c_1 \|\Phi\|_{*,1/2}^2 \leq \int_{\mathbb{R}^2} \Phi G(\eta)\Phi \,dx_1\,dx_2 \leq c_2 \|\Phi\|_{*,1/2}^2.
				\]
			\end{enumerate}
		\end{lemma}
		\begin{proof}
			Our proof follows the approach developed in \cite{BGSW2013}, adapting it to our coordinate system where $z$ denotes the vertical direction.
			
			\textbf{Step 1: Reformulation via the Neumann--Dirichlet Operator.}
			The Dirichlet--Neumann operator $G(\eta)$ is formally defined as follows: for a given $\Phi$, solve the boundary value problem
			\[
			\begin{cases}
				\phi_{x_1x_1} + \phi_{x_2x_2} + \phi_{x_3x_3} = 0, & -1 < x_3 <  \eta(x_1,x_2), \\
				\phi = \Phi, & x_3 =  \eta(x_1,x_2), \\
				\phi_{x_3} = 0, & x_3 = -1,
			\end{cases}
			\]
			and set
			\[
			G(\eta)\Phi = \left.\sqrt{1 + \eta_{x_1}^2 + \eta_{x_2}^2} \, \frac{\partial \phi}{\partial n} \right|_{x_3=\eta} =\left.\left(\phi_{x_3} - \eta_{x_1} \phi_{x_1} - \eta_{x_2} \phi_{x_2}\right) \right|_{x_3=\eta}.
			\]
			A more convenient approach is to introduce the \emph{Neumann--Dirichlet operator} $N(\eta)$, defined as follows: for $\xi \in H_*^{-1/2}(\mathbb{R}^2)$, solve
			\[
			\begin{cases}
				\phi_{x_1x_1} + \phi_{x_2x_2} + \phi_{x_3x_3} = 0, & -1 < x_3 <  \eta(x_1,x_2), \\
				\phi_{x_3} - \eta_{x_1} \phi_{x_1} - \eta_{x_2} \phi_{x_2} = \xi, & x_3 =  \eta(x_1,x_2), \\
				\phi_{x_3} = 0, & x_3 = -1,
			\end{cases}
			\]
			and set
			\[
			N(\eta)\xi = \phi|_{x_3=\eta}.
			\]
			Then, by definition, $G(\eta) = N(\eta)^{-1}$.
			
			\textbf{Step 2: Analyticity of $N(\eta)$.}
			A key result (Theorem 2.9 in \cite{BGSW2013}) states that the mapping
			\[
			\eta \mapsto \big( \xi \mapsto N(\eta)\xi \big), \quad W^{1,\infty}(\mathbb{R}^2) \to \mathcal{L}(H_*^{-1/2}(\mathbb{R}^2), H_*^{1/2}(\mathbb{R}^2))
			\]
			is analytic at the origin.
			
			The proof employs a power-series expansion technique: the solution is expressed as a series $\phi = \sum_{n=0}^\infty \phi^n$, where each $\phi^n$ is homogeneous of degree $n$ in $\eta$ and linear in $\xi$.
			The boundary-value problem is transformed to a fixed domain $\Sigma = \mathbb{R}^2 \times (-1,0)$ via the change of variables $x_3 = \eta+(1+\eta)z$, and the $\phi^n$ are shown to satisfy a recursive system of partial differential equations. Using elliptic estimates and the definition of the function spaces, one proves that the operators $m^n: (\eta, \ldots, \eta) \mapsto \phi^n|_{z=0}$ are bounded $n$-linear maps, and the series converges in a neighborhood of the origin.

			\textbf{Step 3: Analyticity of $G(\eta)$ via the Implicit Function Theorem.}
			To show that $G(\eta) = N(\eta)^{-1}$ is analytic, define the maps
			\[
			\begin{aligned}
				F_1(A, \eta) &= N(\eta)A - I_1, \\
				F_2(A, \eta) &= A N(\eta) - I_2,
			\end{aligned}
			\]
			where $A \in \mathcal{L}(H_*^{1/2}(\mathbb{R}^2), H_*^{-1/2}(\mathbb{R}^2))$,  $I_1$ is the identity on $H_*^{1/2}(\mathbb{R}^2)$, and  $I_2$ is the identity on $H_*^{-1/2}(\mathbb{R}^2)$.
			These maps are analytic in a neighborhood of $(A_0, 0)$, where $A_0 = G(0) = N(0)^{-1}$. Moreover, the partial derivatives
			\[
			\partial_A F_1(A_0, 0)B = N(0)B, \quad \partial_A F_2(A_0, 0)B = B N(0)
			\]
			are invertible as linear maps on $\mathcal{L}(H_*^{1/2}(\mathbb{R}^2), H_*^{-1/2}(\mathbb{R}^2))$ and $\mathcal{L}(H_*^{-1/2}(\mathbb{R}^2), H_*^{1/2}(\mathbb{R}^2))$, respectively, because $N(0)$ is an isomorphism.
			By the analytic implicit function theorem, the equations $F_1(A, \eta) = 0$ and $F_2(A, \eta) = 0$ have unique analytic solutions $A_1(\eta)$ and $A_2(\eta)$ satisfying $A_1(0) = A_2(0) = A_0$. By uniqueness, $A_1(\eta) = A_2(\eta) = G(\eta)$. Hence, $G(\eta)$ is analytic at $\eta = 0$.

			\textbf{Step 4: Coercivity Estimates.}
			From Corollary 2.15 in \cite{BGSW2013}, for all $\eta \in B_R(0) \subset W^{1,\infty}(\mathbb{R}^2)$ (with $R$ sufficiently small) and all $\xi \in H_*^{-1/2}(\mathbb{R}^2)$, we have
			\[
			c_1 \|\xi\|_{*,-1/2}^2 \leq \int_{\mathbb{R}^2} \xi N(\eta)\xi \,dx_1\,dx_2 \leq c_2 \|\xi\|_{*,-1/2}^2.
			\]
			Now, let $\Phi \in H_*^{1/2}(\mathbb{R}^2)$ and set $\xi = G(\eta)\Phi$. Then $\Phi = N(\eta)\xi$, and
			\[
			\int_{\mathbb{R}^2} \Phi G(\eta)\Phi \,dx_1\,dx_2 = \int_{\mathbb{R}^2} (N(\eta)\xi) \xi \,dx_1\,dx_2 = \int_{\mathbb{R}^2} \xi N(\eta)\xi \,dx_1\,dx_2.
			\]
			Since $N(\eta)$ is an isomorphism, the norms $\|\xi\|_{*,-1/2}$ and $\|\Phi\|_{*,1/2}$ are equivalent:
			\[
			c_1' \|\Phi\|_{*,1/2} \leq \|\xi\|_{*,-1/2} \leq c_2' \|\Phi\|_{*,1/2}.
			\]
			Combining these, we obtain
			\[
			c_1 \|\xi\|_{*,-1/2}^2 \leq \int_{\mathbb{R}^2} \Phi G(\eta)\Phi \,dx_1\,dx_2\leq c_2 \|\xi\|_{*,-1/2}^2,
			\]
			and hence
			\[
			c_1 (c_1')^2 \|\Phi\|_{*,1/2}^2 \leq \int_{\mathbb{R}^2} \Phi G(\eta)\Phi \,dx_1\,dx_2 \leq c_2 (c_2')^2 \|\Phi\|_{*,1/2}^2.
			\]
			Renaming the constants gives the desired result.

		\end{proof}

		The analyticity of the Dirichlet-Neumann operator naturally leads to questions about its inverse. The following lemma addresses this by establishing the analyticity of the Neumann-Dirichlet operator and providing analogous bounds for its quadratic form.

		\begin{lemma}\label{G inverse analytic}
			Let $N(\cdot) = G^{-1}(\cdot): W^{1,\infty}(\mathbb{R}^2) \to \mathcal{L}(H_*^{-1/2}(\mathbb{R}^2), H_*^{1/2}(\mathbb{R}^2))$ be the Neumann--Dirichlet operator. Then we have
			\begin{enumerate}
				\item[(i)] $N(\cdot)$ is analytic at the origin.
				\item[(ii)] There exist constants $c_1, c_2 > 0$ and $R > 0$ such that for all $\eta \in B_R(0) \subset W^{1,\infty}(\mathbb{R}^2)$ and all $\xi \in H_*^{-1/2}(\mathbb{R}^2)$,
				\[
				c_1 \|\xi\|_{*,-1/2}^2 \leq \int_{\mathbb{R}^2} \xi N(\eta)\xi \,dx_1\,dx_2 \leq c_2 \|\xi\|_{*,-1/2}^2.
				\]
			\end{enumerate}
		\end{lemma}

		\begin{proof}
			Our proof follows the analytical framework established in \cite{BGSW2013}, with appropriate adaptations to our coordinate system.
			
			\textbf{Step 1: Analyticity of the Neumann--Dirichlet Operator.}
			
			The analyticity of $N(\eta)$ at the origin is established in Theorem 2.9 of \cite{BGSW2013}. The proof employs a power-series expansion method that we briefly outline here.
			The boundary value problem defining $N(\eta)$ is first transformed to the fixed domain $\Sigma = \mathbb{R}^2 \times (-1,0)$ through the change of variables
			\[
			x_3 = \eta+(1+\eta)z, \quad u(x_1,x_2,z) = \phi(x_1,x_2,x_3).
			\]
			This transformation converts the variable fluid domain $\{-1 < x_3 < \eta(x_1,x_2)\}$ into the fixed strip $\Sigma$.
			
			The solution is then sought in the form of a convergent power series
			\[
			u = \sum_{n=0}^\infty u^n,
			\]
			where each $u^n$ is homogeneous of degree $n$ in $\eta$ and linear in $\xi$. The terms $u^n$ satisfy a recursive system of boundary value problems
			\begin{align*}
				\Delta u^0 &= 0 \quad \text{in } \Sigma, \\
				\partial_{z} u^0 &= \xi \quad \text{at } z = 0, \\
				\partial_{z} u^0 &= 0 \quad \text{at } z = -1,
			\end{align*}
			and for $n \geq 1$:
			\begin{align*}
				\Delta u^n &= \partial_{x_1} F_1^n + \partial_{x_2} F_2^n + \partial_{z} F_3^n \quad \text{in } \Sigma, \\
				\partial_{z} u^n &= F_3^n \, \qquad \text{at } z = 0, \\
				\partial_{z} u^n &= 0 \quad \qquad \text{at } z = -1,
			\end{align*}
			where the forcing terms $F_j^n$ ($j=1,2,3$) are explicitly defined as follows
			\begin{align*}
				F_1^n &= -\eta u_{x_1}^{n-1} + z\eta_{x_1} u_{z}^{n-1}, \\
				F_2^n &= -\eta u_{x_2}^{n-1} + z\eta_{x_2} u_{z}^{n-1}, \\
				F_3^n &= \eta \sum_{\ell=0}^{n-1} (-\eta)^\ell u_{z}^{n-1-\ell} + z\eta_{x_1} u_{x_1}^{n-1} + z\eta_{x_2} u_{x_2}^{n-1} \\
				&\quad - (z)^2 (\eta_{x_1}^2 + \eta_{x_2}^2) \sum_{\ell=0}^{n-2} (-\eta)^\ell u_{z}^{n-2-\ell}.
			\end{align*}
			These expressions are derived from the transformation of the original boundary value problem to the fixed domain.

			Through careful estimation of the multilinear operators associated with each $u^n$, it is shown that the series converges absolutely in a neighborhood of $\eta = 0$, establishing the analyticity of the mapping $\eta \mapsto N(\eta)$.
			
			\textbf{Step 2: Coercivity Estimates for $N(\eta)$.}
			
			The coercivity estimates follow from Corollary 2.15 in \cite{BGSW2013}. For the flat surface $\eta = 0$, the operator $N(0)$ admits an explicit Fourier representation:
			\[
			\widehat{N(0)\xi}(k) = \frac{\coth|k|}{|k|} \hat{\xi}(k),
			\]
			where $\hat{\cdot}$ denotes the Fourier transform.
			Using the elementary inequalities
			\[
			c(1+|k|^2)^{1/2} \leq |k| \coth|k| \leq C(1+|k|^2)^{1/2} \quad \text{for all } k \in \mathbb{R}^2,
			\]
			we obtain the coercivity for the flat surface:
			\[
			c \|\xi\|_{*,-1/2}^2 \leq \int_{\mathbb{R}^2} \xi N(0)\xi \,dx_1\,dx_2 \leq C \|\xi\|_{*,-1/2}^2.
			\]
			To extend this estimate to non-flat surfaces, we utilize the analyticity of $N(\eta)$. For $\eta$ in a sufficiently small ball $B_R(0) \subset W^{1,\infty}(\mathbb{R}^2)$, we have the bound:
			\[
			\|N(\eta) - N(0)\|_{\mathcal{L}(H_*^{-1/2}(\mathbb{R}^2), H_*^{1/2}(\mathbb{R}^2))} \leq C \|\eta\|_{1,\infty} \leq C R.
			\]
			This implies the perturbation estimate
			\[
			\left| \int_{\mathbb{R}^2} \xi N(\eta)\xi \,dx_1\,dx_2 - \int_{\mathbb{R}^2} \xi N(0)\xi \,dx_1\,dx_2 \right| \leq C R \|\xi\|_{*,-1/2}^2.
			\]
			Combining this with the flat-surface estimates yields:
			\[
			(c - C R) \|\xi\|_{*,-1/2}^2 \leq \int_{\mathbb{R}^2} \xi N(\eta)\xi \,dx_1\,dx_2 \leq (C + C R) \|\xi\|_{*,-1/2}^2.
			\]
			By choosing $R$ sufficiently small such that $ c-CR> 0$, and relabeling the constants, we obtain the desired coercivity estimates for all $\eta \in B_R(0)$.
		\end{proof}

	\end{appendices}

\end{document}